\documentclass{article}
\usepackage{arxiv}

\usepackage[utf8]{inputenc}
\usepackage[T1]{fontenc}

\usepackage{amsmath,amssymb,amsfonts,amsthm}
\usepackage{bm}
\DeclareMathOperator*{\argmin}{argmin}

\usepackage{graphicx}
\usepackage{subcaption}
\graphicspath{{./Figures/}}

\usepackage{algorithm}
\usepackage{algorithmic}

\usepackage{booktabs}
\usepackage{nicefrac}
\usepackage{multicol}
\usepackage{pifont}
\usepackage{url}
\usepackage{hyperref}
\usepackage{doi}
\usepackage{microtype}
\usepackage{textcomp}
\usepackage{color,xcolor}
\usepackage{setspace}
\usepackage{comment}
\usepackage{tikz}
\usetikzlibrary{calc,arrows}

\newcommand{\vertiii}[1]{{\left\vert\kern-0.25ex\left\vert\kern-0.25ex\left\vert #1
    \right\vert\kern-0.25ex\right\vert\kern-0.25ex\right\vert}}

\newtheorem{theorem}{Theorem}

\newtheorem{lemma}{Lemma}
\newtheorem{assumption}{Assumption}
\newtheorem{proposition}{Proposition}
\newtheorem{remark}{Remark}

\usepackage{todonotes}

\def\grad{\nabla}

\def\cG{\mathcal{G}}

\def\cN{\mathcal{N}}

\def\cV{\mathcal{V}}

\def\smskip{\smallskip}

\def\texitem#1{\par\smskip\noindent\hangindent 25pt
               \hbox to 25pt {\hss #1 ~}\ignorespaces}

\newcommand{\BEAS}{\begin{eqnarray*}}
\newcommand{\EEAS}{\end{eqnarray*}}
\newcommand{\BEA}{\begin{eqnarray}}
\newcommand{\EEA}{\end{eqnarray}}
\newcommand{\BEQ}{\begin{eqnarray}}
\newcommand{\EEQ}{\end{eqnarray}}
\newcommand{\BIT}{\begin{itemize}}
\newcommand{\EIT}{\end{itemize}}
\newcommand{\BNUM}{\begin{enumerate}}
\newcommand{\ENUM}{\end{enumerate}}

\newcommand{\BA}{\begin{array}}
\newcommand{\EA}{\end{array}}

\newcommand{\reals}{\mathbb{R}}

\newif\ifpagenumbering
\pagenumberingtrue

\pagenumberingfalse

\newsavebox{\theorembox}
\newsavebox{\lemmabox}
\newsavebox{\defnbox}
\newsavebox{\corollarybox}
\newsavebox{\assbox}
\savebox{\theorembox}{\noindent\bf Theorem}
\savebox{\lemmabox}{\noindent\bf Lemma}
\savebox{\defnbox}{\noindent\bf Definition}
\savebox{\corollarybox}{\noindent\bf Corollary}
\newtheorem{defn}{\usebox{\defnbox}}

\def\BibTeX{{\rm B\kern-.05em{\sc i\kern-.025em b}\kern-.08em
    T\kern-.1667em\lower.7ex\hbox{E}\kern-.125emX}}

\title{AdGT: Decentralized Gradient Tracking with Adaptive Per-Agent
Stepsizes%
\thanks{This work was partially supported by the Research Council of Finland (Grant 354523). We acknowledge the computational resources provided by the Aalto Science-IT project.}%
\thanks{$\dagger$This work has been submitted to the IEEE for possible publication. Copyright may be transferred without notice, after which this version may no longer be accessible.}
}

\author{
Diyako Ghaderyan\\
Department of Information and Communications Engineering, Aalto University\\
02150 Espoo, Finland\\
\texttt{diyako.ghaderyan@aalto.fi}
\And
Stefan Werner\\
Department of Electronic Systems, Faculty of Information Technology and Electrical Engineering,\\
Norwegian University of Science and Technology (NTNU), 7034 Trondheim, Norway\\
and Department of Information and Communications Engineering, Aalto University, 02150 Espoo, Finland\\
\texttt{stefan.werner@ntnu.no}
}

\date{}

\renewcommand{\shorttitle}{AdGT: Decentralized Gradient Tracking with Adaptive Per-Agent
Stepsizes}

\begin{document}

\maketitle
{
\begin{abstract}
In decentralized optimization, gradient-tracking
methods typically rely on a single global stepsize. This choice
can be conservative when agents have local objectives with
different smoothness constants, since the stepsize must remain
stable for the agent with the largest smoothness constant. This
paper proposes AdGT, a decentralized gradient-tracking method
in which each agent adapts its own stepsize using local gradient
variation and a single global safety factor. The method reduces
fixed-stepsize tuning effort and allows agents to exploit local
smoothness information during the iterations. For smooth and
strongly convex local objectives over undirected networks, we
prove that the analyzed AdGT update converges linearly to
the exact consensus optimizer.  We also study {two adaptive stepsize updates that use changes in} the gradient-tracking direction{. We characterize when the corresponding candidate determines the stepsize and prove conditional lower and upper stepsize bounds and linear convergence under an additional relative tracking-disagreement condition}.
 Experiments on logistic regression, ridge regression, synthetic quadratic problems, and a linear-regression benchmark against state-of-the-art decentralized solvers show that AdGT often {reaches a given accuracy in fewer iterations or gradient evaluations} than tuned fixed-stepsize {GT and the tested baselines}, especially under heterogeneous local smoothness. {In the topology experiments, each tested AdGT update uses one common safety factor across all graphs, whereas the fixed GT stepsize is tuned separately for each graph and seed}.
\end{abstract}}

\keywords{
Adaptive stepsize,  decentralized optimization, gradient tracking, consensus.
}

\section{Introduction}
Recent   advances in artificial intelligence and communication technologies have accelerated the adoption of distributed optimization techniques to solve large-scale problems. These methods allow individual agents to perform local computations and share information with their neighbors, which is particularly advantageous when data is stored across distributed locations or partitioned across multiple servers to improve training efficiency. Key challenges in decentralized optimization include maintaining data privacy, ensuring algorithmic scalability, and achieving robustness to network and system-level disturbances. Distributed optimization has found widespread applications in areas such as machine learning~\cite{lian2017can, tang2018d} and control systems~\cite{nedic2018distributed}.  Large-scale prediction problems in {finance provide further examples} of data-intensive learning applications~\cite{jin2025high}. A common objective in distributed optimization is to solve the consensus problem, which can be formulated as:
\begin{align}
\label{eqz1}
{{\bm{x}^{*}}\in}\argmin\limits_{{\bm{x}}{\in\reals^p}}{{\bar{f}}}({\bm{x}}){\triangleq} \frac{1}{n} \sum\limits_{i=1}^{n} f_{i}({\bm{x}}),
\end{align}
where the objective function { $\bar{f}$ represents  the 
average of } all individual cost functions $\{f_i\}_{i=1}^n$, with $f_{i} : \mathbb{R}^{p}\rightarrow \mathbb{R}$ being the private function of agent $i$.

Building on the foundational work of~\cite{tsitsiklis1986distributed}, the authors of~\cite{nedic2009distributed} introduced a distributed (sub)gradient method for solving the consensus problem in~\eqref{eqz1}. When each local function $f_i$ is convex, the method achieves sublinear convergence. This relatively slow rate is a known limitation of subgradient methods and is largely due to the use of diminishing stepsizes. To improve convergence,~\cite{jakovetic2014fast} proposed an accelerated method that assumes bounded and Lipschitz continuous gradients. By incorporating Nesterov’s acceleration technique, they obtained a convergence rate of $\mathcal{O}(\log(k)/k^2)$, which significantly outperforms standard subgradient methods. Under similar assumptions, EXTRA~\cite{shi2015extra} was developed for smooth convex functions with Lipschitz gradients. This method uses a fixed stepsize and introduces a correction term involving the difference between consecutive gradients. It achieves a convergence rate of $\mathcal{O}(1/k)$, which improves to linear convergence when the objective functions are strongly convex.  {In addition, the work in \cite{shi2015proximal} investigated a proximal-gradient method for decentralized optimization.} In~\cite{yuan2018exact}, the authors proposed the \emph{Exact Diffusion} method that was shown to achieve linear convergence under standard assumptions.

In addition to gradient-based methods, several distributed algorithms have been developed based on the \emph{alternating direction method of multipliers (ADMM)}. These approaches are known for their strong convergence guarantees and ability to handle more complex and structured optimization problems. Notable examples include~\cite{wei20131,shi2014linear,mokhtari2016dqm,mokhtari2016decentralized}, which address general consensus optimization settings. Furthermore, several ADMM-based methods have been proposed for solving {composite convex objectives}~\cite{aybat2015asynchronous,aybat2017distributed}, as well as for problems involving complicated constraints~\cite{aybat16,aybat2019distributed}.

Gradient tracking methods (GT)~\cite{di2015distributed,di2016next,nedic2017achieving,qu2017harnessing,xu2015augmented} form a widely studied class of algorithms in decentralized optimization. These methods introduce an auxiliary variable to estimate the global gradient by locally tracking the average gradient across the network. This mechanism enables their performance to closely match that of centralized algorithms, which often exhibit linear convergence rates. These techniques have been effectively adapted to a wide range of challenging settings, including directed and time-varying communication graphs, asynchronous computation, composite optimization problems, and even nonconvex objectives; see, for instance,~\cite{nedic2017achieving,xi2017add,xu2017convergence,pu2021distributed,koloskova2021improved,sun2022distributed,tian2020achieving,xin2018linear,lu2021optimal,liu2024decentralized,jiang2025convergence,wu2025effectiveness,sun2023decentralized,ye2025generalization}. More recently, \cite{ghaderyan2023fast} proposed a new approach that incorporates gradient tracking in an implicit way by introducing a novel momentum term, specifically designed for use over directed networks. To improve convergence rates and computational complexity, a variety of accelerated decentralized gradient-based methods have been proposed; see, for example,~\cite{qu2019accelerated, kovalev2021lower, rogozin2019optimal, ye2023multi, ye2020decentralized, li2020decentralized, xu2020accelerated, li2020revisiting, li2024accelerated, scaman2017optimal}. 
 
However, a single fixed stepsize {can be overly conservative when the local
objectives have different smoothness constants or condition numbers. In
standard fixed-stepsize GT analyses, the admissible common stepsize depends on
global objective and network constants and, in particular, must be small enough
for the agent with the largest local smoothness constant}. {Consequently, agents
with smoother local objectives cannot take advantage of larger steps, and a
small number of restrictive agents can slow the convergence of the entire
network. Gradient-tracking methods with} uncoordinated, agent-dependent
stepsizes {have been studied} over both undirected and directed
graphs~{\cite{xu2015augmented,nedic2017geometrically,xin2019frost}.
Nevertheless, selecting and tuning a separate fixed stepsize for every agent
remains difficult and becomes impractical in large decentralized networks.
These observations} motivate the following question:

\begin{center}
\itshape
{Can} each agent automatically {adapt its own} stepsize at every iteration {using
only local information, without line search or manual per-agent tuning, while
preserving convergence and improving practical performance}?
\end{center}

{In centralized optimization, this} question has been {studied for convex and
strongly convex objectives. In particular,~\cite{malitsky2020adaptive}
proposed an adaptive gradient method based on two main principles: the stepsize
should not increase too quickly, and it should not be too large for the
observed local smoothness. The method estimates local smoothness from
successive iterates and gradients, without requiring function-value
evaluations or line search. Because the observed local smoothness can be
smaller than a global smoothness bound, the method can use larger steps while
retaining convergence guarantees}. Several extensions have {developed less
conservative stepsize updates or adaptive proximal-gradient methods; see}, for
example,
{\cite{latafat2024adaptive,malitsky2024adaptive,zhou2025adabb}}.

Adaptive stepsize {methods} have also received increasing {attention} in federated
learning, including client adaptivity and adaptive server
optimizers~\cite{wang2021local,wang2022communication,reddi2020adaptive}{.
In particular, building on~\cite{malitsky2020adaptive},
\cite{kim2023adaptive} allows each client to automatically} select its own
stepsize {using a local smoothness estimate}.

{These developments suggest that adaptation to local smoothness may be
particularly useful in decentralized optimization, where agents can have
different local objectives and communicate only with their neighbors. This
leads to a second question}:

\begin{center}
\itshape
{Can adaptive per-agent stepsizes} provide especially large practical gains in
decentralized networks {when the local objectives have different smoothness
constants}?
\end{center}

In this paper, we address these questions by {proposing an adaptive gradient-tracking method, termed AdGT, and by studying when per-agent stepsize adaptation} is most beneficial. {During the preparation of this work, a related decentralized adaptive method was proposed in~\cite{chen2024distributed}. Both methods select per-agent stepsizes using locally available information, but their updates differ. The update in~\cite{chen2024distributed} combines a local smoothness estimate, the norm of the gradient tracker, and a fixed growth factor. AdGT uses a history-dependent growth-control candidate and, in \eqref{Stepsize_Adgt_y_i} and \eqref{Stepsize_Adgt_y_i_f_i}, a tracking-based candidate defined from the change in the gradient-tracking direction}. Both analyses establish uniform {lower and upper stepsize bounds}; see Lemma~\ref{lemma:alpha_k_bounded} and \cite[Lemma~1]{chen2024distributed}. The {method in~\cite{chen2024distributed} converges to a neighborhood of the optimizer, whereas the main AdGT update converges linearly to the exact consensus optimizer}. In our analysis, the consensus, {gradient-tracking}, and optimality errors {form one error vector satisfying a linear recursion; contraction of this recursion drives all three errors} to zero.

\textit{Contributions:} {The main contributions are as follows.}
\begin{itemize}

\item We propose \emph{AdGT}, a decentralized gradient-tracking method {in
which each agent adapts its} stepsize using local gradient {information and a
common safety factor. AdGT requires neither line search nor per-agent stepsize
tuning. Relative to standard adapt-then-combine GT, it requires no additional
communication round and exchanges the same two $p$-dimensional vectors per
iteration. For its main update, we establish uniform lower and upper stepsize
bounds and} prove exact linear convergence to the consensus optimizer {under the
stated assumptions}.

\item We {develop two additional stepsize updates based on changes in the
gradient-tracking direction. We characterize when the tracking-based candidate
determines the next stepsize and derive a one-step bound showing how the
observed change in this direction limits the selected stepsize. Under an
additional relative tracking-disagreement condition, we establish uniform
stepsize bounds and prove linear convergence for both updates.}

\item {We evaluate AdGT on real and synthetic problems under different graph
topologies and levels of local smoothness heterogeneity. Empirically, AdGT
reduces graph-dependent tuning and often reaches a given residual level in
fewer iterations} than tuned fixed-stepsize GT{. The quadratic experiments show
larger gains when only a small subset of agents has large local smoothness
constants: the performance gap between AdGT and GT increases, whereas the
corresponding AdGD--GD gap changes much less.}

\end{itemize}

The {remainder of this} paper is organized as follows. Section~II {presents AdGT} and its adaptive per-agent stepsize {updates}. Section~III {gives the convergence analysis, including the conditional results for \eqref{Stepsize_Adgt_y_i} and \eqref{Stepsize_Adgt_y_i_f_i}}. Section~IV {reports} numerical results on real and synthetic problems {under different network topologies and levels of heterogeneity}. Section~V concludes the paper and discusses future {work.}

\textit{Notation:}
  Throughout the paper, we use the notation $\|\cdot\|$ to denote the Euclidean $2$-norm for vectors and the Frobenius norm for matrices. In particular, for any matrices $A$ and $B$, \[
\|AB\| \leq \vertiii{A}\,\|B\|,
\]
where $\|\cdot\|$ denotes the Frobenius norm and $\vertiii{\cdot}$ the spectral  norm.\footnote{See \cite{horn2012matrix}, Problem~5.6.P20 and the proof of Theorem~7.4.10.1.} 

\section{Design of AdGT}
%Our goal is to solve 
Consider the consensus optimization problem described by \eqref{eqz1} over a communication network represented by an \emph{undirected} connected graph $\mathcal{G} = (\mathcal{V}, \mathcal{E})$. The graph consists of a set of nodes (agents) $\mathcal{V} = \{1, 2, \ldots, n\}$ and a set of edges $\mathcal{E}$. Each node $i \in \mathcal{V}$ has an individual cost function $f_i : \mathbb{R}^p \rightarrow \mathbb{R}$, which is only accessible to that node. 

For each node, $i \in \mathcal{V}$, $\mathcal{N}_i = \{j \in \mathcal{V} : (j, i) \in \mathcal{E}\} \cup \{i\}$ is the set of neighbors of agent $i$. i.e., all nodes that can communicate directly with node $i$. Additionally, $W = [w_{ij}] \in \mathbb{R}^{n \times n}$ is a symmetric doubly-stochastic matrix satisfying:
\begin{align} \label{eqz5}
     w_{ij}=\left\{
                \begin{array}{ll}
              > 0, &\quad  j \in {\cN_{i},}\\
                  0,&  \quad \text{otherwise{;}}
                \end{array}
              \right.
    \quad \sum\limits_{{j\in \cV}} w_{ij}=1, \quad \forall~i{\in\cV}.
\end{align} 
which ensures that the weights sum to one across both rows and columns.

Throughout the paper, we make the following assumptions.
\begin{assumption}\label{Assum::network}
{$\cG$} is undirected and  connected.
\end{assumption}

\begin{remark}\label{Remark: tilde_W_lambda}
Let
{$\lambda_1\geq\lambda_2\geq\cdots\geq\lambda_n$}
denote the eigenvalues of $W$. Since {$\mathcal G$} is connected and $W$ is
{symmetric and doubly stochastic with positive diagonal entries},
$\lambda_1=1$ and $|\lambda_i|<1$ for $i=2,\ldots,n$. Define
\[
 {\widetilde W
 =
 W-\frac{\bm 1\bm 1^\top}{n},
 \qquad
 \lambda
 :=
 \max\{|\lambda_2|,|\lambda_n|\}.}
\]
Then
\[
 {\vertiii{\widetilde W}=\lambda<1.}
\]

\end{remark}
\begin{assumption}\label{Assum:Lsmooth}
{For every $i\in\cV$, the local function $f_{i}$ is $L_i$-smooth, i.e., %such that 
it is differentiable 
with a Lipschitz 
gradient:}
\begin{equation}\label{eqz2}
\Vert \nabla f_{i}({\bm{x}})-\nabla f_{i}({\bm{x}'})\Vert \leq L_i \Vert {\bm{x}} -{\bm{x}'}  \Vert,\quad{\forall~{\bm{x}},\bm{x}' \in \reals^{p}.}
\end{equation}
\end{assumption}

\begin{remark}\label{remark:L_max}
The Lipschitz constant of $\nabla f$ is given by $L = \max_i L_i$, where $L_i$ is the Lipschitz constant of the local function $f_i$.
\end{remark}

\begin{assumption}\label{Assum:Strongly_connected}
{For all $i\in\cV$,} $f_{i}$ is $\mu_i$-strongly convex, i.e.,

\begin{equation}\label{eqz3}
 f_{i}({\bm{y}}) \geq  f_{i}({\bm{x}}) + \nabla f_{i}({\bm{x}})^\top ({\bm{y}}-{\bm{x}}) + \dfrac{\mu_i}{2} \parallel {\bm{y}} -{\bm{x}}  \Vert^2
\end{equation}
for some $\mu_i>0$ and all ${\bm{x}},\bm{y}\in\reals^{p}$. 
\end{assumption}
\begin{remark}
Under Assumption~\ref{Assum:Strongly_connected}, the optimal solution to \eqref{eqz1} is unique, denoted by $\bm{x}^*$.
\end{remark}
\begin{remark}\label{remark:mu_min}
Throughout the analysis, we use
\[
    \mu \triangleq \min_{i\in\mathcal{V}} \mu_i .
\]
\end{remark}
\begin{defn}
\label{def:f}
Define
  {$\bm{\mathrm{x}} \triangleq \left[ \bm{x}_1,\ldots,\bm{x}_n\right]^\top \in \reals^{{n\times p}}$ and 
  $\bm{\mathrm{y}} \triangleq \left[ \bm{y}_1,\ldots,\bm{y}_n\right]^\top \in \reals^{{n\times p}}$,}
where $\bm{x}_i$, $\bm{y}_i \in \reals^{p} $ are the local variables of agent $i\in\cV\triangleq\{1,\ldots,n\}$, and in an algorithmic framework, their values at iteration {$k\geq 0$} are $\bm{x}_i^{k}$ and  $\bm{y}_i^{k}$ for $i\in\cV$.  Let $f:\reals^{{n \times p}}\to\reals$ be a function of local variables $\{\bm{x}_i\}_{i\in\cV}$ such that  $f(\bm{\mathrm{x}})\triangleq \sum_{i\in\cV} f_i(\bm{x}_i)$ for $\bm{\mathrm{x}}\in\reals^{{n\times p} }$ and {$\grad f(\bm{\mathrm{x}}) \triangleq \left[ \grad f_1(\bm{x}_1),...,\grad f_n(\bm{x}_n)\right]^\top \in \reals^{{n \times p}}$,} where $\grad f_{i}(\bm{x}_i) \in \mathbb{R}^{p}$ denotes the gradient of $f_{i}$ at $ \bm{x}_i\in\reals^p$. Define
{$\bm{\bar{x}}\triangleq \tfrac{1}{n}\bm{1}^\top\bm{\mathrm{x}}\in\reals^{1\times p}$, 
$\bm{\bar{y}}\triangleq \tfrac{1}{n}\bm{1}^\top\bm{\mathrm{y}}\in\reals^{1\times p}$, 
$\bm{\bar{\mathrm{x}}}\triangleq \bm{1}\bm{\bar{x}}\in\reals^{n\times p}$, 
$\bm{\bar{\mathrm{y}}}\triangleq \bm{1}\bm{\bar{y}}\in\reals^{n\times p}$, 
$\grad f(\bm{\bar{x}})\triangleq \left[\grad f_1(\bm{\bar{x}}),\ldots,\grad f_n(\bm{\bar{x}})\right]^\top\in\reals^{n\times p}$, 
$\grad f(\bm{x}^\ast)\triangleq \left[\grad f_1(\bm{x}^\ast),\ldots,\grad f_n(\bm{x}^\ast)\right]^\top\in\reals^{n\times p}$}.
\end{defn}

We next propose a decentralized optimization algorithm, AdGT, to solve 
{the consensus optimization problem in~\eqref{eqz1}.}

\begin{algorithm}
\caption{AdGT}
\label{algori1}
\textbf{Input:} {$\bm{x}_{i}^{0}\in\reals^{p},~\forall\,i\in\cV$}
\begin{algorithmic}[1]
\STATE {$\bm{y}_{i}^{0}= \nabla f_{i}(\bm{x}_{i}^{0})$, $\alpha_{i}^{0}> 0$, $\gamma > 0$,
        $\theta_{i}^{0}=0$, for $i\in\cV$}
\FORALL{$k=0,1,\ldots$}
  \STATE {Each agent $i\in\mathcal{V}$ independently performs:}
  \STATE {$\bm{x}_{i}^{k+1} = \sum\limits_{j\in\cN_i} w_{ij}\!\left(\bm{x}_{j}^{k}
          - \alpha_j^{k}\, \bm{y}_{j}^{k} \right)$} \label{AdGT_line_cons}
  \STATE {$\bm{y}_{i}^{k+1} = \sum\limits_{j\in\cN_i} w_{ij}\bm{y}_{j}^{k}
          + \nabla f_{i}(\bm{x}_{i}^{k+1}) - \nabla f_{i}(\bm{x}_{i}^{k})$} \label{AdGT_line_track}
\STATE {Update $\alpha_i^{k+1}$ using \eqref{alpha_nound_L_F}, \eqref{Stepsize_Adgt_y_i}, or \eqref{Stepsize_Adgt_y_i_f_i}.}
  \STATE {$\theta_{i}^{k+1}= \dfrac{\alpha_{i}^{k+1}}{\alpha_{i}^{k}}$}
\ENDFOR
\end{algorithmic}
\end{algorithm}
\vspace*{-1mm}
In Algorithm~\ref{algori1}, $\alpha_i^0>0$ is {the initial local} stepsize used
in the first primal update{. Each agent then updates its stepsize locally} using
one of the {adaptive updates} below. For $k\ge1$,
$\theta_i^k=\alpha_i^k/\alpha_i^{k-1}$ {is the ratio of two consecutive
stepsizes and appears in the term that controls stepsize growth.}

{
\subsection{Adaptive {Stepsize Design}}

{The adaptive stepsize updates} in Algorithm~\ref{algori1} {build on} the
centralized adaptive gradient method {in}~\cite{malitsky2020adaptive}, {which uses}
\begin{align}{\label{Adap_Cent_Step}
 \alpha^{k+1}
 =
 \min\left\{
 \frac{\|\bm x^{k+1}-\bm x^k\|}
      {2\|\nabla f(\bm x^{k+1})-\nabla f(\bm x^k)\|},
 \sqrt{1+\theta^k}\,\alpha^k
 \right\},}
\end{align}
where $\theta^k=\alpha^k/\alpha^{k-1}$. The {first term selects the stepsize
using the local smoothness observed from two consecutive iterates and
gradients, while the second limits how quickly the stepsize can increase.

In decentralized gradient tracking, the primal update uses the local gradient
tracker $\bm y_i^k$ rather than only the local gradient. This direction depends
on both the local objective and information received from neighboring agents.
Consequently, local gradient variation alone may not fully describe the change
in the direction used in the primal update.

A related decentralized adaptive method was proposed
in~\cite{chen2024distributed} with}
\begin{align}{\label{Adap_Dist_DM}
 \alpha_i^{k+1}
 =
 \min\left\{
 \frac{1}{2L_{f_i}^k},
 \frac{1}{\|\bm y_i^{k+1}\|},
 \sqrt{2}\,\alpha_i^k
 \right\},}
\end{align}
where
\begin{align}{\label{L_smooth_L_f_i}
 L_{f_i}^k
 =
 \frac{
 \|\nabla f_i(\bm x_i^{k+1})-\nabla f_i(\bm x_i^k)\|
 }{
 \|\bm x_i^{k+1}-\bm x_i^k\|
 }.}
\end{align}
The {quantity $L_{f_i}^k$ estimates the local Lipschitz constant of
$\nabla f_i$ from two consecutive iterates and gradients. Throughout the
paper, we call it the \emph{local smoothness estimate}.}
Update~\eqref{Adap_Dist_DM} combines this estimate, a term based on the norm of
the gradient tracker, and a fixed limit on stepsize growth. Its analysis
guarantees convergence to a neighborhood of the {optimizer}. In our experiments,
{we observed numerical instability for this method in some settings with few
local} samples or sparse network connectivity.

{Motivated by~\eqref{Adap_Cent_Step}, AdGT introduces a common safety factor and
a history-dependent limit on stepsize growth. The update} used in the {main
convergence analysis is}
\begin{align}{\label{alpha_nound_L_F}
 \alpha_i^{k+1}
 =
 \min\left\{
 \frac{1}{2\gamma L_{f_i}^k},
 \sqrt{1+\theta_i^k}\,\alpha_i^k
 \right\},}
\end{align}
{for every agent and $k\geq0$, with} $\theta_i^0=0$. Here, $\gamma>0$ is {the
common} safety factor; see Remark~\ref{remark:gamma_explain}. {We call
$1/(2\gamma L_{f_i}^k)$ the \emph{local-gradient candidate} and
$\sqrt{1+\theta_i^k}\alpha_i^k$ the \emph{growth-control candidate}}. The first
uses the local smoothness estimate, while the second {prevents a rapid increase
in the stepsize}. At $k=0${,} update~\eqref{alpha_nound_L_F} {becomes}
\begin{align}
 {\alpha_i^1
 =
 \min\left\{
 \frac{1}{2\gamma L_{f_i}^0},
 \alpha_i^0
 \right\}.
 \label{eq:alpha_initial}}
\end{align}

{Update~\eqref{alpha_nound_L_F} uses the change in the local gradient. To also
account for the change in the direction used in the primal update, define}
\begin{align*}
 {L_{y_i}^k
 =
 \frac{\|\bm y_i^{k+1}-\bm y_i^k\|}
      {\|\bm x_i^{k+1}-\bm x_i^k\|}.}
\end{align*}
{
We call $L_{y_i}^k$ the \emph{tracking-direction variation estimate} and
$1/(2\gamma L_{y_i}^k)$ the \emph{tracking-based candidate}. Using this
candidate together with the growth-control candidate gives
} 
\begin{align}{\label{Stepsize_Adgt_y_i}
 \alpha_i^{k+1}
 =
 \min\left\{
 \frac{1}{2\gamma L_{y_i}^k},
 \sqrt{1+\theta_i^k}\,\alpha_i^k
 \right\}.}
\end{align}
{The gradient trackers preserve the average-gradient identity}
\begin{align}{\label{eq:ybar_identity}
 \bar{\bm y}^k
 =
 \frac{1}{n}\sum_{i=1}^{n}\nabla f_i(\bm x_i^k);}
\end{align}
{see~\cite{qu2017harnessing}. Thus,} $L_{y_i}^k$ measures {the change in the
local gradient-tracking direction relative to the local primal change}.

{Update~\eqref{Stepsize_Adgt_y_i} isolates} the effect of {the tracking-based
candidate. However, the local-gradient change and the consensus correction in
the tracking update may partly cancel. In this case,} $L_{y_i}^k$ {can become
small and its candidate can become large. We therefore also consider}
\begin{align}{\label{Stepsize_Adgt_y_i_f_i}
 \alpha_i^{k+1}
 =
 \min\left\{
 \frac{1}{2\gamma L_{f_i}^k},\;
 \frac{1}{2\gamma L_{y_i}^k},\;
 \sqrt{1+\theta_i^k}\,\alpha_i^k
 \right\},}
\end{align}
{which retains the local-gradient candidate while also using the change in the
gradient-tracking direction.

For any of the three updates, if
$\|\bm x_i^{k+1}-\bm x_i^k\|=0$ for some $k\geq1$, we set}
$\alpha_i^{k+1}=\alpha_i^k$. If {this} denominator is nonzero but {the numerator
defining $L_{f_i}^k$ or $L_{y_i}^k$ is zero}, the corresponding {candidate is
interpreted} as $+\infty$ and {omitted from} the minimum. The {quantities required
by the selected update at} $k=0$ are assumed to be well defined.

{Since updates~\eqref{Stepsize_Adgt_y_i} and
\eqref{Stepsize_Adgt_y_i_f_i} both contain the tracking-based candidate, we
refer to them collectively as the \emph{tracking-aware updates}. The main
convergence theorem covers update~\eqref{alpha_nound_L_F}.
Remark~\ref{rem:tracking_candidate} explains the role of the tracking-based
candidate, while Proposition~\ref{prop:tracking_bounds_convergence} provides
conditional stepsize bounds and convergence results for the tracking-aware
updates.

{
\begin{remark}
\label{rem:tracking_candidate}
For $\|\bm x_i^{k+1}-\bm x_i^k\|>0$, define
\[
 \bm g_i^k
 =
 \nabla f_i(\bm x_i^{k+1})-\nabla f_i(\bm x_i^k),
 \qquad
 \bm r_i^k
 =
 \sum_jw_{ij}(\bm y_j^k-\bm y_i^k),
\]
where $\bm r_i^k$ is the consensus correction produced by averaging the
neighboring tracking variables. The tracking update gives
\[
 \bm y_i^{k+1}-\bm y_i^k
 =
 \bm g_i^k+\bm r_i^k,
\]
and hence
\[
 (L_{y_i}^k)^2-(L_{f_i}^k)^2
 =
 \frac{
 2\langle\bm g_i^k,\bm r_i^k\rangle+\|\bm r_i^k\|^2
 }{
 \|\bm x_i^{k+1}-\bm x_i^k\|^2
 }.
\]
Therefore,
\[
 \frac{1}{2\gamma L_{y_i}^k}
 \le
 \frac{1}{2\gamma L_{f_i}^k}
 \quad\Longleftrightarrow\quad
 2\langle\bm g_i^k,\bm r_i^k\rangle+\|\bm r_i^k\|^2\ge0.
\]

For update~\eqref{Stepsize_Adgt_y_i_f_i}, the tracking-based candidate is
selected, possibly with a tie, when
\[
 2\langle\bm g_i^k,\bm r_i^k\rangle+\|\bm r_i^k\|^2\ge0
 \quad\text{and}\quad
 \frac{1}{2\gamma L_{y_i}^k}
 \le
 \sqrt{1+\theta_i^k}\,\alpha_i^k.
\]
When both inequalities are strict, it uniquely determines the next stepsize.
Moreover, update~\eqref{Stepsize_Adgt_y_i_f_i} contains the two candidates in
\eqref{alpha_nound_L_F} and adds the tracking-based candidate. Thus, when the
two updates are evaluated at the same state with the same
$\alpha_i^k$, $\theta_i^k$, and $\gamma$,
update~\eqref{Stepsize_Adgt_y_i_f_i} cannot produce a larger next stepsize.

Update~\eqref{Stepsize_Adgt_y_i_f_i} also satisfies
\[
 2\gamma\alpha_i^{k+1}
 \|\bm y_i^{k+1}-\bm y_i^k\|
 \le
 \|\bm x_i^{k+1}-\bm x_i^k\|.
\]
Thus, a larger change in the gradient-tracking direction relative to the local
primal change places a smaller upper bound on the next stepsize.

For update~\eqref{Stepsize_Adgt_y_i}, the local-gradient candidate is absent.
Therefore, the tracking-based candidate is selected, possibly with a tie, when
\[
 \frac{1}{2\gamma L_{y_i}^k}
 \le
 \sqrt{1+\theta_i^k}\,\alpha_i^k,
\]
and the same one-step bound holds. Unlike
\eqref{Stepsize_Adgt_y_i_f_i}, however,
\eqref{Stepsize_Adgt_y_i} does not retain the local-gradient candidate when
$L_{y_i}^k$ becomes small. These statements concern one update from a common
state and do not order the complete trajectories or convergence rates of the
different updates.
\end{remark}}}

Section~\ref{Simulation:section_1} reports which candidate is selected,
Remark~\ref{re:lower_clipping} introduces clipping, and
Proposition~\ref{prop:tracking_bounds_convergence} gives conditional stepsize
bounds and convergence results for the original unclipped updates.

\begin{remark}\label{remark:gamma_explain} The parameter $\gamma$ is a {common} safety factor, not {the fixed stepsize used by standard GT. Increasing} $\gamma$ {decreases the local-gradient and tracking-based candidate stepsizes, while the selected stepsizes remain agent- and iteration-dependent. Section~III gives sufficient worst-case conditions on} $\gamma$. In the experiments{,} $\gamma$ is {selected from a small common set, whereas the} fixed GT stepsize {is tuned separately for each graph}. \end{remark}
 For $k\geq0$, the primal and gradient-tracking recursions of AdGT can be
written compactly as
\begin{align}
\bm{\mathrm{x}}^{k+1} &=  W \Big( \bm{\mathrm{x}}^k - D^k \bm{\mathrm{y}}^k \Big) \label{eqComp1} \\
\bm{\mathrm{y}}^{k+1} &= W \bm{\mathrm{y}}^k +  \nabla f(\bm{\mathrm{x}}^{k+1})-\nabla f(\bm{\mathrm{x}}^k)\label{eqComp2}
\end{align}
where $D^k=\operatorname{diag}(\alpha_1^k,\ldots,\alpha_n^k)$ contains the
agents' stepsizes. The matrix $D^0$ contains the {initial} stepsizes used in the first update, whereas
the convergence analysis uses the adaptive stepsizes $D^k$ for $k\geq1$.

{
Relative to standard adapt-then-combine GT, AdGT requires no extra communication
round: at iteration \(k\), agent \(j\) broadcasts the two \(p\)-dimensional
vectors \(\bm{x}_j^k-\alpha_j^k\bm{y}_j^k\) and \(\bm{y}_j^k\). Thus, AdGT does
not increase the order of the vector communication payload. Its additional
local cost consists only of vector differences, norm evaluations, and scalar
updates for \(\alpha_i^k\) and \(\theta_i^k\), which is \(O(p)\) per agent and
iteration; \(L_{y_i}^k\) is computed from quantities already available after the
tracking update.
}

{
\section{Main Results}
 This section proves linear convergence of AdGT with update~\eqref{alpha_nound_L_F}. The proof {first} establishes uniform stepsize bounds and then {combines} the consensus, {gradient-tracking}, and optimality errors {in a linear recursion}. The resulting rate is a conservative worst-case guarantee.

{Proposition~\ref{prop:tracking_bounds_convergence} applies the same error recursions to} updates~\eqref{Stepsize_Adgt_y_i} and {\eqref{Stepsize_Adgt_y_i_f_i} under an additional relative tracking-disagreement condition. This condition is not implied by Assumptions~\ref{Assum::network}--\ref{Assum:Strongly_connected}}.

We first show that the adaptive stepsizes generated by AdGT remain uniformly
bounded, a property that is essential for the convergence analysis.

\begin{lemma}
\label{lemma:alpha_k_bounded}
Suppose Assumptions~\ref{Assum:Lsmooth} and~\ref{Assum:Strongly_connected} hold.
Let $\{\alpha_i^k\}_{k\ge1}$ be generated by Algorithm~\ref{algori1} with update~\eqref{alpha_nound_L_F}. Assume that $L_{f_i}^0$ is well
defined for every agent $i$, and that the {initial} stepsize satisfies
$\alpha_i^0\ge 1/(2\gamma L_i)$. Then, for every $i\in\mathcal V$ and all $k\ge1$,
\begin{equation}\label{bound_alpha}
  \frac{1}{2 \gamma L} \leq \alpha_i^{k} \leq \frac{1}{2 \gamma \mu}.  
\end{equation}
\end{lemma}

The proof is given in Appendix~\ref{app:proof_alpha_k_bounded}. Defining
\[
{\alpha_{\max}\triangleq
\sup_{i\in\mathcal V,\; k\ge1}\alpha_i^k.}
\]
Lemma~\ref{lemma:alpha_k_bounded} gives
\begin{equation}\label{define_bound_alpha_max}
  \frac{1}{2 \gamma L} \leq \alpha_{\max} \leq \frac{1}{2 \gamma \mu}.  
\end{equation}
\begin{remark}
\label{re:lower_clipping}
{For numerical diagnostics, we also report} the lower-clipped {form}
\begin{equation}{\label{Stepsize_Adgt_clipped}
 \alpha_{i,\mathrm{clip}}^{k+1}
 =\max\{\underline{\alpha}_i,\alpha_{i,(12)}^{k+1}\},}
\end{equation}
where {$\alpha_{i,(12)}^{k+1}$ denotes the right-hand side of
\eqref{Stepsize_Adgt_y_i_f_i} and $\underline{\alpha}_i>0$ is a small numerical
floor}. {This} floor is only a numerical safeguard {and is not used in the
conditional result below}.
\end{remark}

For completeness, we restate another technical result; see
\cite[Lemma~10]{qu2017harnessing} for a proof.
\begin{lemma}\label{lemma_max_L_mu}
Let $\alpha \in \bigl(0,\tfrac{2}{L}\bigr)$ and define 
\[
\eta \triangleq \max \left\{ |1-L\alpha|,\; |1-\mu\alpha| \right\}.
\]
If Assumptions~\ref{Assum:Lsmooth} and~\ref{Assum:Strongly_connected} hold, then for all $\bm{x}\in \mathbb{R}^p$,
\begin{equation*}
\bigl\| \bm{x}-\alpha \frac{1}{n}\!\!\sum_{i=1}^n \nabla f_i(\bm{x}) - \bm{x}^{*}\bigr\| 
\;\leq\; \eta \,\|\bm{x}-\bm{x}^{*}\|.
\end{equation*}
\end{lemma}
The next lemma collects the three coupled inequalities used in the convergence
proof; they control, respectively, the consensus error, the gradient-tracking
error, and the optimality error.

\begin{lemma}
\label{lemma:one_step_recursion}
Suppose Assumptions~\ref{Assum::network}--\ref{Assum:Strongly_connected} hold. Let
$\gamma \geq \tfrac{L}{\mu}$. Then, for all $k\geq1$, the AdGT iterates satisfy
{
\begin{align}
    \|\bm{\mathrm{x}}^{k+1}- \bar{\bm{\mathrm{x}}}^{k+1}\| 
    &\leq (\lambda+\lambda \alpha_{\max} L)\|\bm{\mathrm{x}}^k- \bar{\bm{\mathrm{x}}}^k\|  
     + \lambda \alpha_{\max} \|\bm{\mathrm{y}}^k- \bar{\bm{\mathrm{y}}}^k \| 
      + \sqrt{n}\,\lambda \alpha_{\max} L\| \bar{\bm{{x}}}^k - \bm{{x}}^* \|,
     \label{eq:main_consensus_recursion}\\[1mm]
   \| \bm{\mathrm{y}}^{k+1}- \bar{\bm{\mathrm{y}}}^{k+1} \| 
   &\leq (\lambda+L\alpha_{\max}) \|\bm{\mathrm{y}}^k- \bar{\bm{\mathrm{y}}}^k \| 
    +  L \big(\vertiii{W-I}+ L\alpha_{\max}\big)\|\bm{\mathrm{x}}^k- \bar{\bm{\mathrm{x}}}^k\|  
   +\sqrt{n}\, L^2 \alpha_{\max} \| \bar{\bm{x}}^k - \bm{x}^* \|,
   \label{eq:main_tracking_recursion}\\[1mm]
\| \bar{\bm{x}}^{k+1} - \bm{x}^*\|  
&\leq  \left(1- \frac{\mu}{2 \gamma L}\right) \,\|\bar{\bm{x}}^{k} -\bm{x}^{*}\| 
     + L \frac{\alpha_{\max}}{\sqrt{n}}  \|\bar{\bm{\mathrm{x}}}^{k} -\bm{\mathrm{x}}^{k}  \| 
       +\frac{\alpha_{\max}}{\sqrt{n}} \|\bm{\mathrm{y}}^{k} -\bar{\bm{\mathrm{y}}}^{k}\| .
      \label{eq:main_average_recursion}
\end{align}
}
\end{lemma}

The proof is given in Appendix~\ref{app:proof_one_step_recursion}. 

The three
errors are mutually coupled, since each bound in
Lemma~\ref{lemma:one_step_recursion} depends on the other two. This is what
makes a comparison-system argument natural: rather than bounding each error in
isolation, we treat them jointly. Once $\gamma$ is chosen so that the
associated comparison matrix has spectral radius strictly smaller than one, the
coupled errors decrease geometrically, which implies exact linear convergence
to the consensus optimizer.

Combining~\eqref{eq:main_consensus_recursion}--\eqref{eq:main_average_recursion}, we obtain
\begin{align}\label{eqz19}
e_{k+1}\leq  \Upsilon~e_{k},\quad \forall~k\geq1,
\end{align}
where $e_k$ and $\Upsilon$ are defined as
{
\begin{align*}
&e_{k}=
\begin{bmatrix}
\|\bm{\mathrm{x}}^k- \bar{\bm{\mathrm{x}}}^k\| \\
 \|\bm{\mathrm{y}}^k- \bar{\bm{\mathrm{y}}}^k \|  \\
\| \bar{\bm{{x}}}^k - \bm{{x}}^* \|
\end{bmatrix},\\
&\Upsilon=
\begin{bmatrix}
\lambda+\lambda \alpha_{\max} L& \lambda \alpha_{\max} & \sqrt{n}\,\lambda \alpha_{\max}L\\
L \big(\vertiii{W-I}+ L\alpha_{\max}\big) & \lambda+L\alpha_{\max} & \sqrt{n}\, L^2 \alpha_{\max} \\
L \frac{\alpha_{\max}}{\sqrt{n}} & \frac{\alpha_{\max}}{\sqrt{n}} & (1- \frac{\mu}{2 \gamma L})
\end{bmatrix}.
\end{align*}%
}
Thus, the convergence analysis reduces to deriving sufficient conditions under which
$\rho(\Upsilon)<1$.

In the following, we focus on the nontrivial case $\lambda>0$. If $\lambda=0$,
the consensus contribution in the first row of the comparison system vanishes,
and the first term in the definition of $\bar\alpha$ can be omitted.
\begin{theorem}
\label{theorem1}
Let the iterates be generated by Algorithm~\ref{algori1} with update~\eqref{alpha_nound_L_F}. Suppose Assumptions~\ref{Assum::network}--\ref{Assum:Strongly_connected} hold, that $L_{f_i}^0$ is well defined for every agent $i$, and that $\alpha_i^0\ge 1/(2\gamma L_i)$. Let
\[
  \gamma \;>\; \max\!\left\{\frac{L}{\mu},\; \frac{1}{2\mu \,\bar{\alpha}} \right\},
\]
where $\bar{\alpha}>0$ is given by
{
\begin{equation}\label{eq:alpha_bar_def}
   \bar{\alpha}
=
   \min\!\left\{
      \frac{(1-\lambda)\,\delta_1}{\lambda\big(L\delta_1+ 1 +\sqrt{n}\,L\delta_2\big)},\;
      \frac{(1-\lambda)-L\vertiii{W-I}\, \delta_1}{L\big(L\delta_1+ 1+\sqrt{n}\,L \delta_2\big)}
   \right\},
\end{equation}}
s.t.
{
\begin{align}\label{detlta_12}
 0 \;<\; \delta_1 \;<\; \frac{(1-\lambda)}{L\vertiii{W-I}}\,, \qquad
  \delta_2 \;>\; \left(\frac{L}{\sqrt{n}}\delta_1+\frac{1}{\sqrt{n}}\right)\frac{L}{\mu^{2}} .
\end{align}
}
Then the spectral radius satisfies $\rho(\Upsilon)<1$.
\end{theorem}

The proof is given in Appendix~\ref{app:proof_theorem1}.

\begin{remark}\label{re:gamma_conservative}
The bound on $\gamma$ in Theorem~\ref{theorem1} is {sufficient but conservative
and is} not a practical tuning rule. {It uses only uniform stepsize bounds;
Remarks~\ref{remark:gama_bound_fin}
and~\ref{re:beta_1_beta_2_extended} give the derivation}. In the experiments,
$\gamma$ {is chosen from a small common set while each agent still adapts its
stepsize;} fixed-stepsize GT {tunes one step for each graph}.
\end{remark}

\begin{remark}\label{remark:gama_bound_fin}
Theorem~\ref{theorem1} gives a constructive way to guarantee
$\rho(\Upsilon)<1$. First, choose $\delta_1$ and $\delta_2$ satisfying the
feasibility conditions in~\eqref{detlta_12}; these values determine
$\bar{\alpha}$ through~\eqref{eq:alpha_bar_def}. Then any $\gamma$ satisfying
{
\[
    \gamma > \max\!\left\{\frac{L}{\mu},\,\frac{1}{2\mu\,\bar\alpha}\right\}
\]
}
guarantees $\rho(\Upsilon)<1$. One explicit feasible choice is
{
\begin{align}\label{value_delta_12}
       \delta_1 = \frac{1-\lambda}{2L\vertiii{W-I}}, 
    \qquad
    \delta_2 =
    2\left(\frac{L}{\sqrt{n}}\delta_1 + \frac{1}{\sqrt{n}}\right)
    \frac{L}{\mu^{2}} .
\end{align}
}
Here, the factor $1/2$ places $\delta_1$ strictly below its admissible upper
bound in~\eqref{detlta_12}, while the factor $2$ places $\delta_2$ strictly
above its admissible lower bound. Thus, \eqref{value_delta_12} yields an
explicit value of $\bar{\alpha}$ and, through the bound above, an explicit
sufficient lower bound on $\gamma$.
\end{remark}
% =============================================================
% CORRECTED VERSION of Remark re:beta_1_beta_2_extended
% =============================================================
\begin{remark}\label{re:beta_1_beta_2_extended}
Continuing from the feasible choice in~\eqref{value_delta_12}, we now derive an
explicit sufficient stepsize condition. Define
{
\[
T_1=
\frac{(1-\lambda)\delta_1}
{\lambda\bigl(L\delta_1+1+\sqrt n L\delta_2\bigr)},
\qquad
T_2=
\frac{(1-\lambda)-L\vertiii{W-I}\delta_1}
{L\bigl(L\delta_1+1+\sqrt n L\delta_2\bigr)}.
\]}
Then \(\bar\alpha=\min\{T_1,T_2\}\). Hence, to guarantee
\(\alpha_{\max}<\bar\alpha\), it is enough to lower bound \(T_1\) and \(T_2\).

Substituting \eqref{value_delta_12} and using $0<1-\lambda<1$, $\lambda<1$,
$\vertiii{W-I}\le 2$, and $\vertiii{W-I}\ge 1-\lambda$ gives
\[
T_1
\ge
\frac{(1-\lambda)^2}
{5L\!\left(1+\frac{2L^2}{\mu^2}\right)},
\qquad
T_2
\ge
\frac{1-\lambda}
{5L\!\left(1+\frac{2L^2}{\mu^2}\right)}.
\]
Since $0<1-\lambda<1$, the first lower bound is smaller, and therefore
\[
\bar\alpha
\ge
\frac{(1-\lambda)^2}
{5L\!\left(1+\frac{2L^2}{\mu^2}\right)}.
\]
Writing $\kappa=L/\mu$ {gives} \[ {\bar\alpha \ge \frac{(1-\lambda)^2} {5\mu\kappa(1+2\kappa^2)}.} \] {Therefore, a sufficient stepsize condition is} \begin{equation}{\label{eq:AdGT_static_compact_simplified} \alpha_{\max} < \frac{(1-\lambda)^2} {5\mu\kappa(1+2\kappa^2)}.} \end{equation} {Since Lemma~\ref{lemma:alpha_k_bounded} gives $\alpha_{\max}\le1/(2\gamma\mu)$, condition \eqref{eq:AdGT_static_compact_simplified} is ensured by} \[ {\gamma > \max\left\{ \kappa,\, \frac{5\kappa(1+2\kappa^2)} {2(1-\lambda)^2} \right\}.} \] {This is a sufficient} worst-case {condition, not a practical tuning rule}. 
\end{remark}}

\begin{remark}\label{re:compare_stepsize_bounds}
For comparison, representative sufficient fixed-stepsize bounds for the methods
in~\cite{qu2017harnessing,nedic2017achieving}, expressed in our notation, are
\[
 \alpha_{\mathrm{QL}}
 =
 \frac{(1-\lambda)^2}{36\mu\kappa^2},
 \qquad
 \alpha_{\mathrm{DIGing}}
 \le
 \frac{(1-\lambda)^2}
 {2\mu\kappa(1+4\sqrt n\,\sqrt\kappa)}.
\]
The sufficient AdGT condition
in~\eqref{eq:AdGT_static_compact_simplified} has the same
$(1-\lambda)^2$ dependence and scales as
$\mathcal{O}((1-\lambda)^2/(\mu\kappa^3))$. Its dependence on $\kappa$ is more
conservative because our analysis uses uniform lower and upper bounds on the
adaptive stepsizes. Since these bounds come from different worst-case analyses,
they do not rank the practical performance of the methods.
\end{remark}

\begin{theorem}
\label{theorem2}
Let the iterates be generated by Algorithm~\ref{algori1} with update~\eqref{alpha_nound_L_F}. Suppose Assumptions~\ref{Assum::network}--\ref{Assum:Strongly_connected} hold, that $L_{f_i}^0$ is well defined for every agent $i$, and that $\alpha_i^0\ge 1/(2\gamma L_i)$. Suppose further that
\[
  \gamma \;>\; \max\!\left\{\frac{L}{\mu},\; \frac{1}{2\mu \,\bar{\alpha}} \right\},
\]
where $\bar{\alpha}>0$ is defined in~\eqref{eq:alpha_bar_def}. Then the sequence $\{\bm{\mathrm{x}}^{k}\}$ converges linearly to
$\bm{\mathrm{x}}^{*}=\bm{1}\bm{x}^*$ after the initialization step, with a rate
arbitrarily close to $\rho(\Upsilon)$.
\end{theorem}

The proof is given in Appendix~\ref{app:proof_theorem2}.

{ 
\begin{proposition}
\label{prop:tracking_bounds_convergence}
Suppose Assumptions~\ref{Assum::network}--\ref{Assum:Strongly_connected} hold
and, for each agent $i$, there is a constant $\xi_i\ge0$ such that
\[
 \left\|
 \sum_jw_{ij}(\bm y_j^k-\bm y_i^k)
 \right\|
 \le
 \xi_i\|\bm x_i^{k+1}-\bm x_i^k\|
\]
whenever $\|\bm x_i^{k+1}-\bm x_i^k\|>0$. We refer to this inequality as the
\emph{relative tracking-disagreement condition}. It requires the norm of the
consensus correction to be at most $\xi_i$ times the change in agent~$i$'s
local decision variable. Assume that the quantities required
at $k=0$ are well defined and
$\alpha_i^0\ge1/(2\gamma L_i)$.

For update~\eqref{Stepsize_Adgt_y_i_f_i}, the generated stepsizes satisfy
\[
 \frac{1}{2\gamma(L_i+\xi_i)}
 \le
 \alpha_i^k
 \le
 \frac{1}{2\gamma\mu_i},
 \qquad
 i\in\mathcal V,\quad k\ge1.
\]
Define
\[
 \underline\alpha
 =
 \frac{1}{2\gamma\max_i(L_i+\xi_i)}.
\]
Form a comparison matrix from $\Upsilon$ in~\eqref{eqz19} by replacing
$\alpha_{\max}$ with $1/(2\gamma\mu)$ and its $(3,3)$ entry with
$1-\mu\underline\alpha$. If
\[
 \gamma\ge\frac{L}{\mu}
\]
and this matrix has spectral radius below one, then the iterates generated by
\eqref{Stepsize_Adgt_y_i_f_i} converge linearly to
$\bm 1\bm x^*$, with a rate arbitrarily close to that spectral radius.

For update~\eqref{Stepsize_Adgt_y_i}, the same lower bound holds. If, in
addition, $\xi_i<\mu_i$ for every agent, then
\[
 \frac{1}{2\gamma(L_i+\xi_i)}
 \le
 \alpha_i^k
 \le
 \frac{1}{2\gamma(\mu_i-\xi_i)},
 \qquad
 i\in\mathcal V,\quad k\ge1.
\]
The same convergence conclusion holds after replacing $\alpha_{\max}$ in
\eqref{eqz19} with
\[
 \frac{1}{2\gamma\min_i(\mu_i-\xi_i)},
\]
keeping the $(3,3)$ entry equal to $1-\mu\underline\alpha$, and requiring
\[
 \gamma
 \ge
 \frac{L}{\min_i(\mu_i-\xi_i)}.
\]
\end{proposition}
The displayed conditions on $\gamma$ ensure the required upper-stepsize bound;
linear convergence also requires the corresponding spectral-radius condition.
The proof is in Appendix~\ref{app:proof_tracking_bounds}. The relative
tracking-disagreement condition must hold whenever the local primal change is
nonzero and is not implied by the standing assumptions.}

\begin{remark} The adaptive stepsize mechanism is not limited to strongly convex objectives, but the present analysis focuses on this case to prove exact linear convergence. Here Assumption~\ref{Assum:Strongly_connected} is a \emph{local} strong-convexity assumption: each $f_i$ is assumed $\mu_i$-strongly convex. This local assumption {gives $L_{f_i}^k\ge\mu_i$, which} is used to obtain the
uniform upper bound on the generated stepsizes in
Lemma~\ref{lemma:alpha_k_bounded}. Strong convexity also gives the contraction factor \(1-\mu\underline{\alpha}\), where \(\underline{\alpha}>0\) is a uniform lower stepsize bound. For merely convex objectives, this linear-rate argument no longer applies directly; a different descent or function-value analysis would be needed, typically yielding a sublinear rate. Extending the analysis to the case where only $\bar f$ is strongly convex is left for future work.
\end{remark}
\section{Numerical Results}
{

This section evaluates whether adaptive per-agent stepsizes reduce tuning effort
and improve convergence relative to carefully tuned fixed-stepsize GT. Unless otherwise stated, AdGT refers to {update~\eqref{alpha_nound_L_F}, which is covered by Theorem~\ref{theorem2}. The convergence guarantees for updates~\eqref{Stepsize_Adgt_y_i} and \eqref{Stepsize_Adgt_y_i_f_i} are conditional on the relative tracking-disagreement condition in Proposition~\ref{prop:tracking_bounds_convergence}}.
We primarily compare with fixed-stepsize Gradient Tracking (GT) and the adaptive
decentralized method of~\cite{chen2024distributed}, referred to as Method-DM.

All methods use the same graph construction, initialization, data split, and stopping criterion within each experiment. The topology-comparison and $\gamma$-sensitivity experiments report the mean residual over seeds $1,2,3,4$. The diagnostic plots for individual stepsizes and active-candidate behavior use one representative realization, seed $2$. GT is tuned separately for each graph and seed, and Method-DM is implemented in its original form from~\cite{chen2024distributed}.

The experiments are organized as follows. Section~\ref{Simulation:section_1}
studies logistic regression across several graph topologies, including
sensitivity to $\gamma$ and active-candidate behavior. 
Section~\ref{Simulation:section_sota} compares AdGT with state-of-the-art
decentralized methods on a linear-regression benchmark.
Section~\ref{Simulation:section_2} studies heterogeneous local smoothness, and
Section~\ref{Simulation:section_3} examines the effect of the number of agents
and network connectivity.
}
{

\subsection{Adaptive stepsizes across network topologies}
\label{Simulation:section_1}

In this subsection, we study the behavior of the proposed adaptive stepsizes on
decentralized logistic regression. This experiment has three goals. First, we
compare AdGT with carefully tuned fixed-stepsize GT over several graph
topologies. Second, we study the sensitivity of AdGT to the safety factor
$\gamma$. Third, we examine which candidate terms are active in the combined
adaptive update~\eqref{Stepsize_Adgt_y_i_f_i}.

Throughout the experiments, AdGT(8){,} AdGT(10){,} and AdGT(12) denote {stepsize updates~\eqref{alpha_nound_L_F}, \eqref{Stepsize_Adgt_y_i}, and \eqref{Stepsize_Adgt_y_i_f_i}, respectively}. The curve labeled clipped AdGT(17) denotes the numerically safeguarded update~\eqref{Stepsize_Adgt_clipped}.

The decentralized logistic regression problem is formulated as
\begin{align}
\label{eq_Logistic_reg}
\bm{x}^\ast \in \argmin_{\bm{x}\in\mathbb{R}^{p}} f(\bm{x})
= \sum_{i=1}^n f_i(\bm{x}),
\end{align}
where each local function is
\begin{align*}
f_i(\bm{x})
=
\sum_{j=1}^{m_i}
\log\left(1+\exp(-y_i^j M_i^j \bm{x})\right)
+
\frac{\rho}{2}\|\bm{x}\|^2 .
\end{align*}
We use the w8a dataset~\cite{chang2011libsvm}, which contains $49,749$
examples with $300$ features. The features are converted to full format,
standardized using Z-score normalization, and augmented with an intercept term,
so $p=301$. The network consists of $n=16$ agents, and each agent receives
$m_i=25$ samples. The regularization parameter is $\rho=10^{-2}$. For seed $s$,
the data split is generated using data seed $100000+s$, and each graph is
generated using its corresponding graph seed.

We test six communication topologies: Star, Cycle, Line, Ladder, and two random
graphs. The two random graphs are generated using the \texttt{Ggen}
code,\footnote{The \texttt{Ggen} code is written by Dr. W. Shi; see
\url{https://sites.google.com/view/wilburshi/home/research/software-a-z-order/graph-tools/ggen}.}
with connectivity parameters $\mathrm{per}=0.2$ for Random1 and
$\mathrm{per}=0.6$ for Random2. The convergence behavior is measured using the
normalized residual
\begin{align*}
r(k)
\triangleq
\frac{\|\bm{\mathrm{x}}^k-\bm{\mathrm{x}}^\ast\|}
{\|\bm{\mathrm{x}}^0-\bm{\mathrm{x}}^\ast\|}.
\end{align*}

\subsubsection{Comparison across graph topologies}

We first compare AdGT with fixed-stepsize GT and Method-DM over six graph
topologies. The GT stepsize is carefully tuned for each graph and seed, and the
best resulting curve is reported.  For AdGT, we report AdGT(8) and AdGT(12){, corresponding to updates~\eqref{alpha_nound_L_F} and \eqref{Stepsize_Adgt_y_i_f_i}, respectively}. The topology figures use $\gamma=15$ for AdGT(8) and $\gamma=4$ for
AdGT(12), with $\alpha_i^0=10^{-3}$ and $\theta_i^0=0$. {These are stable
operating points for the two updates. Remark~\ref{rem:tracking_candidate} gives
the candidate-level explanation: when $L_{y_i}^k/L_{f_i}^k$ is sufficiently
large, the tracking-based candidate with a smaller value of $\gamma$ can still be no
larger than the local-gradient candidate with a larger $\gamma$. The common
$\gamma=4$ comparison in the active-candidate study below isolates this effect
under the same safety factor}. For fixed-stepsize GT, we perform a fine grid search separately for
each graph and seed and report the best stable run. The selected stepsizes are
listed in Table~\ref{tab:gt_stepsize_centers}, and their variation across graphs
and seeds illustrates the graph- and data-dependent tuning required by
fixed-stepsize GT, which the per-agent adaptive stepsizes of AdGT avoid.

\begin{figure}[!ht]
    \centering
    \subfloat[Star]{%
        \includegraphics[width=0.32\textwidth]{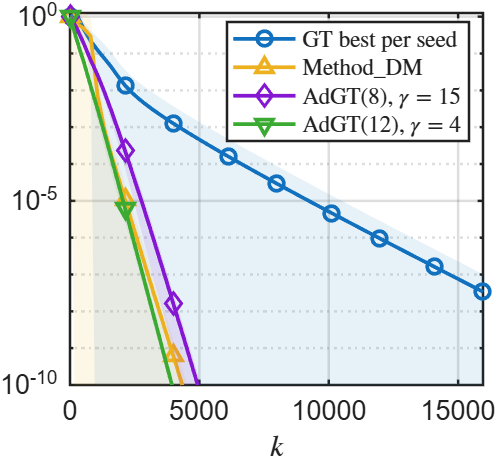}%
    }
    \hfill
    \subfloat[Cycle]{%
        \includegraphics[width=0.32\textwidth]{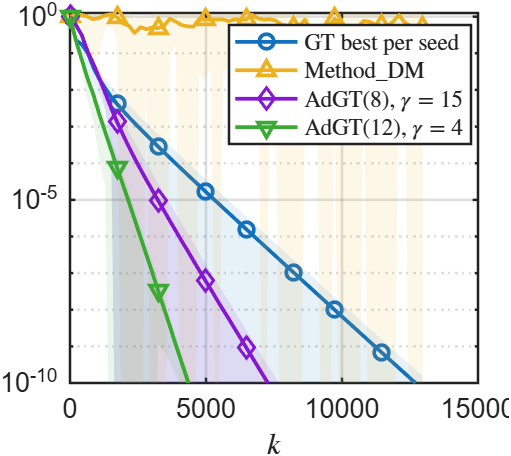}%
    }
   % \vspace{1mm}
     \hfill
    \subfloat[Line]{%
        \includegraphics[width=0.32\textwidth]{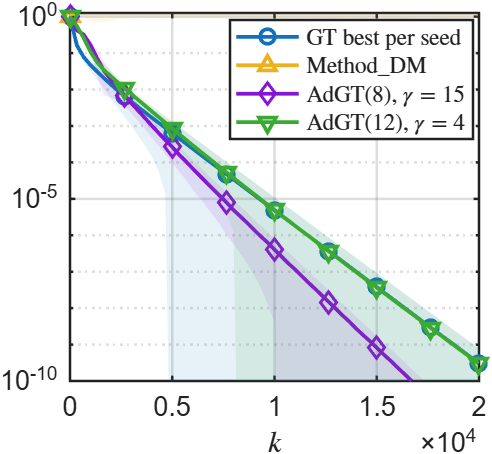}%
    }
    
    \hfill
    \subfloat[Ladder]{%
        \includegraphics[width=0.32\textwidth]{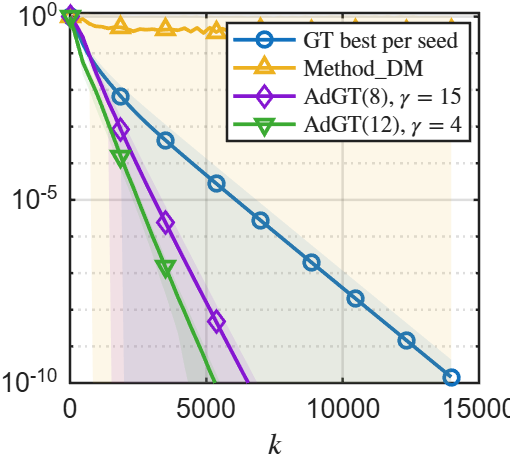}%
    }
    %\vspace{1mm}
 \hfill
    \subfloat[Random1, $\mathrm{per}=0.2$]{%
        \includegraphics[width=0.32\textwidth]{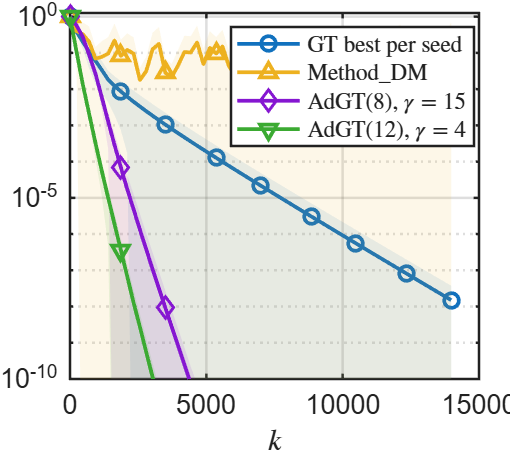}%
    }
    \hfill
    \subfloat[Random2, $\mathrm{per}=0.6$]{%
        \includegraphics[width=0.32\textwidth]{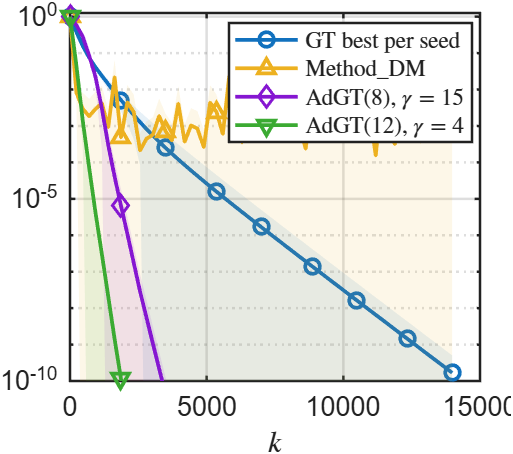}%
    }
    \caption{Mean residual over seeds $1,2,3,4$ for decentralized logistic regression on six network topologies. {AdGT(8) and AdGT(12) each use one common value of $\gamma$ across all six topologies, whereas the fixed} GT {stepsize} is tuned separately for each graph and seed.}
    \label{Logistic_fig}
\end{figure}

Figure~\ref{Logistic_fig} shows that AdGT matches or outperforms carefully
tuned fixed-stepsize GT on the tested topologies, where the GT baseline uses the
favorable per-graph, per-seed stepsizes of Table~\ref{tab:gt_stepsize_centers}.
Method-DM does not reach the same residual level within the plotted iterations
on Cycle, Ladder, Random1, and Random2, while the AdGT variants decrease the
residual more rapidly. The gains are most visible on these harder topologies.
On the Line graph, AdGT(8) {reaches the same residual level in fewer iterations} than tuned GT, whereas AdGT(12) gives comparable performance.

\begin{table}[ht]
\centering
\caption{Best stable fixed GT stepsizes found by fine grid search for each graph and seed.}
\label{tab:gt_stepsize_centers}
\begin{tabular}{c|cccccc}
\hline
Seed & Star & Cycle & Line & Ladder & Random1 & Random2 \\
\hline
1 & 0.16 & 0.32 & 0.13 & 0.26 & 0.21 & 0.21 \\
2 & 0.08 & 0.15 & 0.11 & 0.14 & 0.11 & 0.13 \\
3 & 0.38 & 0.16 & 0.09 & 0.19 & 0.10 & 0.22 \\
4 & 0.28 & 0.37 & 0.10 & 0.40 & 0.26 & 0.21 \\
\hline
\end{tabular}
\end{table}

Two points are important. First, AdGT(8){, which is} covered by the main convergence {theorem, already} matches or outperforms tuned GT {on} the tested topologies{. Thus, the observed} gain does not {depend on adding the tracking-based candidate}. Second, AdGT(12) {includes} the tracking-based candidate $1/(2\gamma L_{y_i}^k)$ {in addition to the candidates in} AdGT(8){. This term can reduce the next stepsize} when the local {smoothness estimate does not capture a large change in the gradient-tracking direction. The} active-candidate study below {examines this behavior}.

\begin{figure}[ht]
\centering
\includegraphics[width=0.45\textwidth]{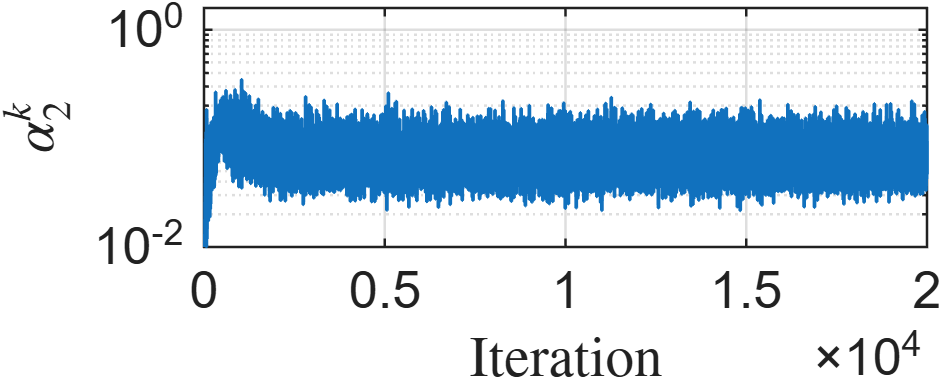}
\hfill
\includegraphics[width=0.45\textwidth]{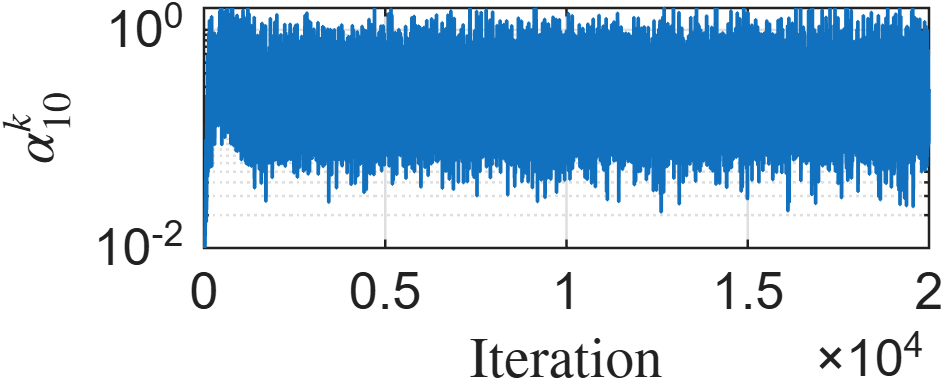}
\caption{Adaptive stepsizes of AdGT(12) for agents $2$ and $10$ on the Line graph, using seed $2$ and $\gamma=4$.}
\label{Stepsise_agent}
\end{figure}

Figure~\ref{Stepsise_agent} shows the stepsize sequences of two representative
agents in the Line graph for the same diagnostic setting used later in the
active-candidate study: seed $2$ and $\gamma=4$. The stepsizes vary across
agents and iterations, which illustrates the main mechanism of AdGT: agents
whose local objective variation permits larger stepsizes can take larger
updates, while agents with more restrictive local behavior use smaller
stepsizes. The persistent variation also shows that AdGT does not simply settle
on one common fixed stepsize.

\subsubsection{Sensitivity to the safety factor $\gamma$}

We next study the sensitivity of AdGT to the safety factor $\gamma$. This
experiment is motivated by the fact that $\gamma$ appears explicitly in the
adaptive stepsize updates. Unlike the fixed stepsize used in GT, however,
$\gamma$ {scales the local-gradient and tracking-based candidate stepsizes but is not itself a stepsize}. The actual stepsizes remain agent-dependent and
iteration-dependent.

The sensitivity experiment was run on all six topologies; Figure~\ref{fig:gamma_sensitivity} reports four representative cases: Cycle, Line, Random2 with $\mathrm{per}=0.6$, and Star. The omitted Ladder and Random1
plots show the same qualitative behavior. The figure compares AdGT(8),
AdGT(10), AdGT(12), and clipped AdGT(17) with $\gamma=15$ and $\gamma=25$. The
curves show the mean residual over seeds $1,2,3,4$. The GT baseline is again
tuned separately for each graph and seed and reported as the best-per-seed
fixed-stepsize baseline. Different line styles and markers are used to
distinguish nearly overlapping adaptive curves.

\begin{figure}[!ht]
\centering

\subfloat[Cycle graph]{%
    \includegraphics[width=0.24\textwidth]{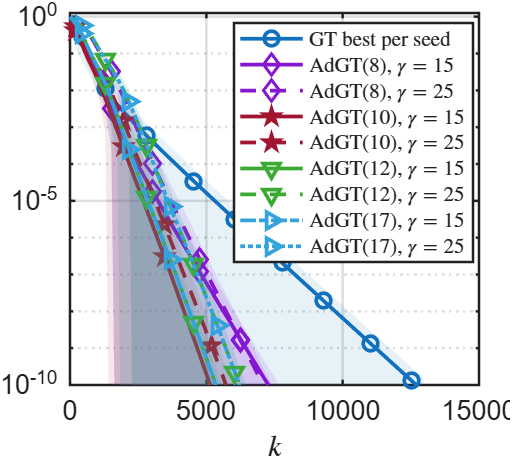}%
}
\hfill
\subfloat[Line graph]{%
    \includegraphics[width=0.23\textwidth]{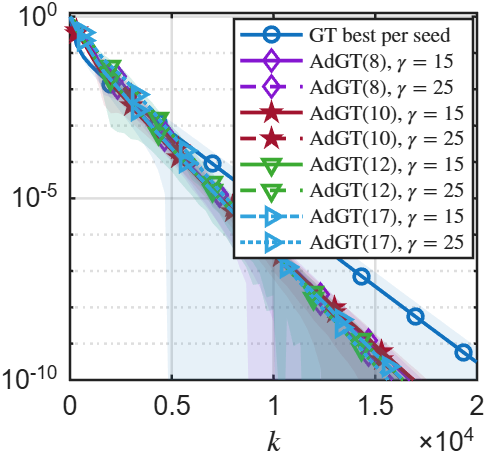}%
}
\hfill
\subfloat[Random2, $\mathrm{per}=0.6$]{%
    \includegraphics[width=0.24\textwidth]{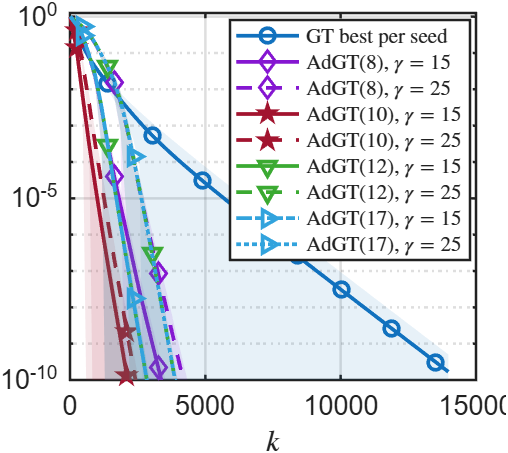}%
}
\hfill
\subfloat[Star graph]{%
    \includegraphics[width=0.23\textwidth]{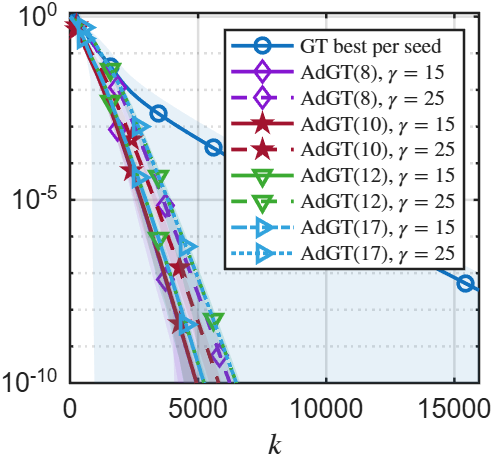}%
}

\caption{Sensitivity of AdGT to $\gamma$ on four representative topologies. Curves show the mean over seeds $1,2,3,4$.}
\label{fig:gamma_sensitivity}
\end{figure}

The results show that, once $\gamma$ is in a stable range, performance
is not highly sensitive to its precise value. The two values $\gamma=15$ and $\gamma=25$ give stable behavior
for all tested adaptive variants and all six topologies, and the four
representative panels in Figure~\ref{fig:gamma_sensitivity} show that these
variants outperform tuned fixed-stepsize GT. As expected, $\gamma=15$ is usually slightly faster and $\gamma=25$ more conservative, but the performance gap between these two values is modest compared with the gap between AdGT and tuned GT. In this stable range, clipped AdGT(17) nearly
overlaps with AdGT(12), which indicates that the numerical floor does not affect
the observed behavior in these sensitivity experiments.

This behavior differs from fixed-stepsize tuning. In GT, the chosen stepsize
sets the update magnitude for every agent and iteration and must be retuned for
each graph and seed. In AdGT, $\gamma$ only scales locally computed candidates;
the actual stepsizes still adapt across agents and iterations through
$L_{f_i}^k$, $L_{y_i}^k$, and the growth-control candidate
$\sqrt{1+\theta_i^k}\alpha_i^k$. Thus, $\gamma$ still needs to be selected, but
it acts as a safety margin rather than as the stepsize itself. Selecting it
from a predefined set of candidate values is less delicate than searching for a
fixed GT stepsize, while still giving better numerical performance in these
experiments. Because $\gamma=15$ and $\gamma=25$ already make the candidates
conservative, this study probes the stable range. The next experiment uses the smaller safety factor $\gamma=4$ to examine a more aggressive regime in which {AdGT(8)} can be unstable {and to study the selection of} the tracking-based candidate {in AdGT(12)}.

\subsubsection{{Candidate selection in AdGT(12)}}

Finally, we examine which candidate {determines the stepsize in} AdGT(12). {In particular, we record how often} the tracking-based {candidate $1/(2\gamma L_{y_i}^k)$ is selected}. Unlike the preceding {sensitivity} study, which uses $\gamma=15$ and $\gamma=25$, {this diagnostic uses the smaller value} $\gamma=4$ to examine a more aggressive {setting. We report} one representative realization with seed{~}$2$, matching the stepsize plots in Figure~\ref{Stepsise_agent}. 

Figure~\ref{fig:active_terms_cycle} reports the {results} on the Cycle graph. {Panel~(a)} compares AdGT(8), AdGT(10), AdGT(12), and clipped AdGT(17) {using the same value} $\gamma=4$. The remaining panels show which candidate determines the AdGT(12) {stepsize. The three candidate stepsizes are} \[ {\sqrt{1+\theta_i^k}\alpha_i^k,\qquad \frac{1}{2\gamma L_{f_i}^k},\qquad \frac{1}{2\gamma L_{y_i}^k},} \] {which are, respectively, the growth-control candidate,} the local-gradient candidate, and the {tracking-based} candidate.

\begin{figure}[!ht]
\centering
\subfloat[Residual comparison]{%
    \includegraphics[width=0.25\textwidth]{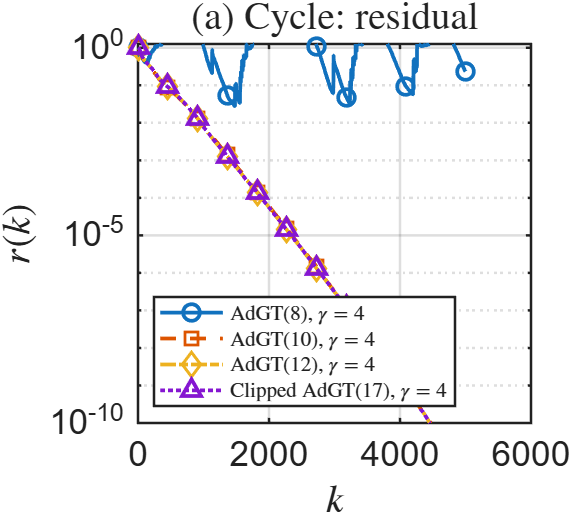}%
}
\hfill
\subfloat[Active-term fractions]{%
    \includegraphics[width=0.25\textwidth]{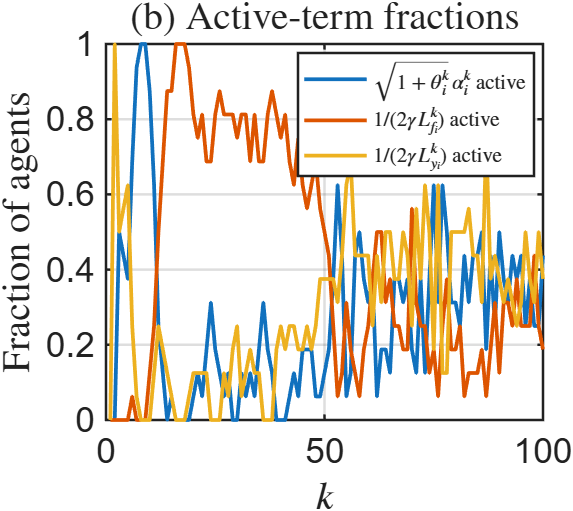}%
}
\hfill
\subfloat[Active term for agent $1$]{%
    \includegraphics[width=0.25\textwidth]{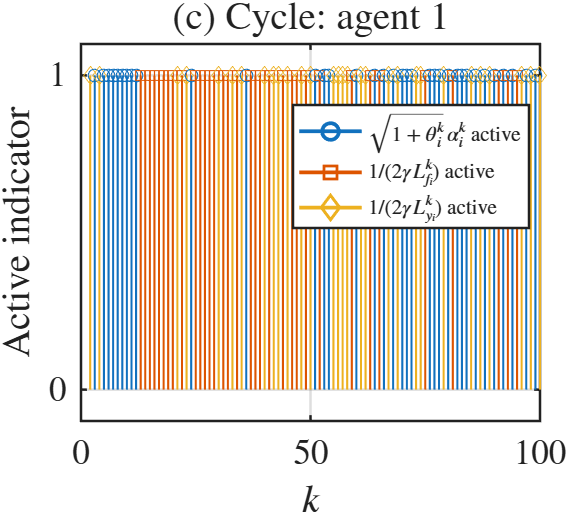}%
}
\hfill
\subfloat[Active term for agent $16$]{%
    \includegraphics[width=0.25\textwidth]{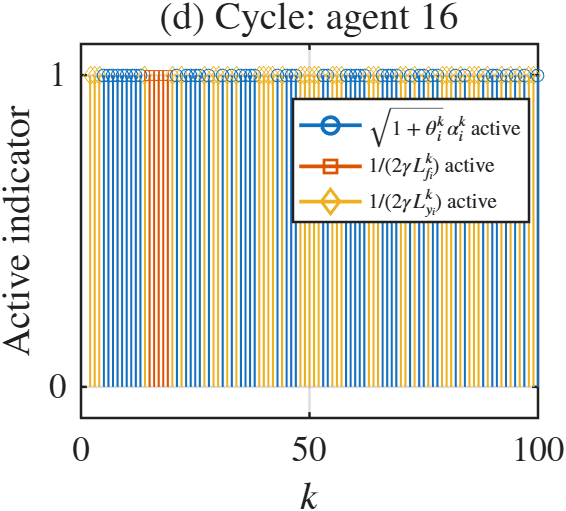}%
}
\caption{Active-candidate behavior on the Cycle graph, using seed $2$ and $\gamma=4$.}
\label{fig:active_terms_cycle}
\end{figure}

Figure~\ref{fig:active_terms_cycle}(a) shows that AdGT(8) exhibits pronounced oscillations and repeated residual increases on the Cycle graph {at} $\gamma=4$. In contrast, AdGT(10), AdGT(12), and clipped AdGT(17) decrease smoothly and reach much smaller residuals under the same safety factor. {In this representative run, the local smoothness estimate alone does not capture all changes in the gradient-tracking direction.}

Figure~\ref{fig:active_terms_cycle}(b) shows that the {tracking-based candidate is selected for} a nonzero fraction of agents{, while} Figures~\ref{fig:active_terms_cycle}(c)--(d) show different {selection} patterns for agents{~}$1$ and{~}$16${. Thus, the tracking-based candidate is not redundant and is selected by different agents at different iterations. These observations support the one-step mechanism in Remark~\ref{rem:tracking_candidate}, but they do not establish a general ordering of the complete trajectories}.

Together, Figures~\ref{Logistic_fig}--\ref{fig:active_terms_cycle} show that
AdGT reduces the need for fixed-stepsize tuning, remains stable under moderate
changes of the safety factor within the tested range, and benefits from the
tracking-based candidate in challenging network topologies.
}

\begin{figure}[!ht]
    \centering
    \includegraphics[width=.65\textwidth]{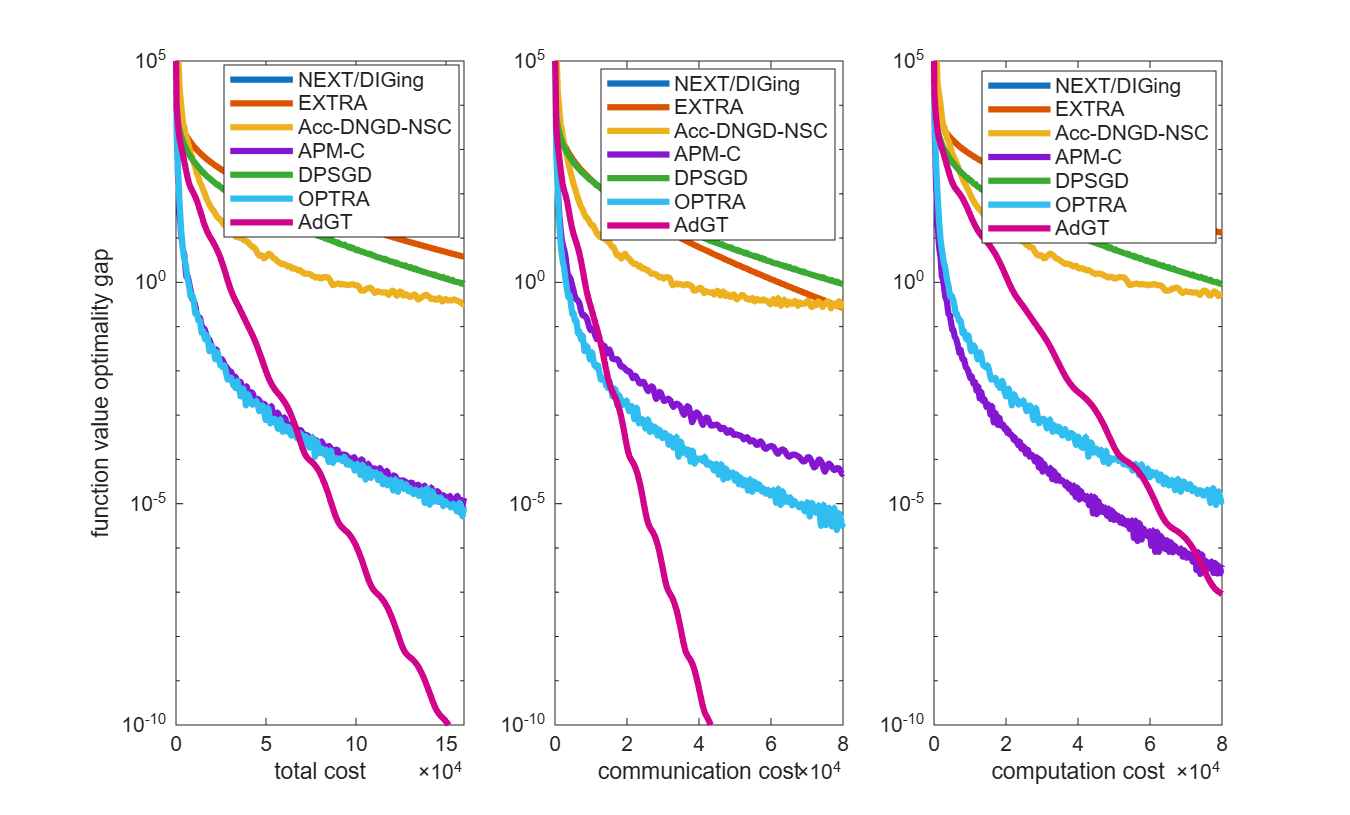}
    \caption{Comparison of decentralized optimization algorithms with AdGT for solving the decentralized linear regression problem.}
    \label{Adgt_OPTRA}
\end{figure}

\subsection{Comparison with State-of-the-Art Decentralized Methods}
\label{Simulation:section_sota}

This section compares the proposed AdGT algorithm with several state-of-the-art
decentralized optimization methods. To ensure a fair and reproducible evaluation, we
follow the experimental framework and parameter settings of Section~5
of~\cite{xu2020accelerated}. All baseline
implementations and data-generation utilities are taken from the official OPTRA
repository.\footnote{\url{https://github.com/YeTian-93/OPTRA}} AdGT is incorporated as an
additional solver while keeping the original data model, network topology, and
hyperparameters unchanged; the number of iterations is moderately increased to better
illustrate long-term convergence behavior. 

The test problem is decentralized linear regression,
\begin{equation}
\min_{\bm{x} \in \mathbb{R}^p} \|A \bm{x} - \bm{b}\|^2,
\end{equation}
where $A = [A_1; A_2; \dots; A_n] \in \mathbb{R}^{nm \times p}$ with
$A_i \in \mathbb{R}^{m \times p}$, and each agent holds its local data
$(A_i,b_i)$. Following the experimental protocol of~\cite{xu2020accelerated},
we use the same problem dimensions and parameter settings as in that work, and
generate $A$ using the correlated Gaussian model of~\cite{agarwal2010fast} with
correlation parameter $\omega=0.95$. The network consists of $20$ agents
connected via an Erd\H{o}s--R\'enyi graph with edge probability $0.1$. The
global smoothness constant $L_f$ is computed from the full dataset as in the
OPTRA protocol and is used by the comparison methods that require it. AdGT does not use $L_f$ to compute its adaptive stepsizes{; its local-gradient
candidate uses a local smoothness estimate computed from consecutive iterates
and gradients}. This benchmark is included to follow the established OPTRA evaluation protocol~\cite{xu2020accelerated} and to compare AdGT with existing decentralized solvers under identical experimental conditions. 

This least-squares benchmark is generally convex rather than strongly convex,
and is used to position AdGT against widely used decentralized solvers under a
standard setup. We compare AdGT with both accelerated and non-accelerated
decentralized algorithms, including OPTRA~\cite{xu2020accelerated},
APM-C~\cite{li2018sharp},
Acc-DNGD-NSC~\cite{qu2019accelerated}, EXTRA~\cite{shi2015extra},
NEXT/DIGing~\cite{di2016next,nedic2017achieving}, and
DPSGD~\cite{lian2017can}. All methods are evaluated under identical conditions
following the protocol of~\cite{xu2020accelerated}. For AdGT, a safety factor
of $\gamma=1$ already yields stable and competitive convergence, while
$\gamma=2$ provides slightly faster empirical performance; hence, we report
results with $\gamma=2$. We observe nearly identical performance for AdGT with
both adaptive stepsize updates~\eqref{alpha_nound_L_F}
and~\eqref{Stepsize_Adgt_y_i_f_i} in this benchmark.  Therefore, we report AdGT with the {local-gradient
update}~\eqref{alpha_nound_L_F} in Figure~\ref{Adgt_OPTRA}.

Figure~\ref{Adgt_OPTRA} {reports the objective} gap versus total computational effort, communication rounds, and gradient evaluations. {Over the plotted range, AdGT reaches a given objective gap with lower total computational effort and with fewer communication rounds and gradient evaluations than the tested} non-accelerated {and accelerated baselines. This is an empirical comparison under the common OPTRA protocol and does not imply a general worst-case complexity ordering}. One possible explanation is that {several baselines use global problem constants or one stepsize for the full run}, whereas AdGT adapts each stepsize to the {local gradient variation observed} along the trajectory. 

\subsection{Adaptive Stepsizes {under Heterogeneous Local Smoothness}}
\label{Simulation:section_2}

{
In this section, we study per-agent {stepsize adaptation} when the local
objectives have different smoothness constants {and condition numbers. A common
GT stepsize must remain stable} for the most restrictive agent{, whereas AdGT
allows} each agent to adapt its own stepsize. {The quadratic model lets us
control the local eigenvalue ranges and inspect the resulting} stepsizes.
Specifically, each agent
\( i \in \{1,\ldots,n\} \) has a local quadratic cost function of the form
}
{
\[
f_i(\bm{x}) := \frac{1}{2} \bm{x}^\top A_i \bm{x} + \bm{b}_i^\top \bm{x},
\]
}
where \( A_i \in \mathbb{S}_{++}^p \) is a positive definite diagonal matrix,
and \( \bm{b}_i \in \mathbb{R}^p \) is a vector. The overall objective is to
minimize the sum of all local functions:
\begin{align}{\label{Objective_quadratic} f(\bm x) = \sum_{i=1}^n \left( \frac12\bm x^\top A_i\bm x+\bm b_i^\top\bm x \right).} \end{align}

{For these quadratic functions,} $L_{f_i}^k$ {measures the local gradient change
relative to} the recent displacement, whereas $L_{y_i}^k$ also {reflects
consensus and gradient-tracking effects.}

To model differences in smoothness across the local objectives, we generate matrices \( A_i \) with varying condition numbers, which directly affect the Lipschitz constants {$L_i=\|A_i\|_2=\lambda_{\max}(A_i)$}. Following the setup in~\cite{mokhtari2016network}, we construct each diagonal matrix \( A_i \) by sampling its entries as follows: the first half of the diagonal (i.e., the first \( p/2 \) entries) is drawn uniformly at random from the set \( \{1, 10^{-1}, \ldots, 10^{-\tau}\} \), and the second half from \( \{1, 10^{1}, \ldots, 10^{\tau}\} \). {Thus, changing $\tau$ varies both local smoothness and conditioning}. The linear terms \( \bm{b}_i^\top \bm{x} \) are added to ensure that the local objective functions have different minimizers. The vectors \( \bm{b}_i \) are generated uniformly at random from the box \( [0, 1]^p \).

We consider a network of \( n = 100 \) nodes, where the dimension of the local variable \( \bm{x} \) is \( p = 20 \). We evaluate the performance under four different scenarios based on the choice of the condition number parameter \( \tau \):
\begin{itemize}
    \item[i)]\label{agent_n_100} All agents are assigned a condition number parameter of \( \tau = 3 \);
    \item[ii)] Half of the agents (i.e., 50 nodes) have \( \tau = 3 \), while the remaining half have \( \tau = 1 \);
    \item[iii)] Ten percent of the agents (i.e., 10 nodes) are assigned \( \tau = 3 \), and the remaining 90 nodes have \( \tau = 1 \);
    \item[iv)] Only 3 agents are assigned \( \tau = 3 \), with all others using \( \tau = 1 \).
\end{itemize}

\begin{figure}[!ht]
    \centering
    \subfloat[Scenario i]{%
        \includegraphics[width=0.24\textwidth]{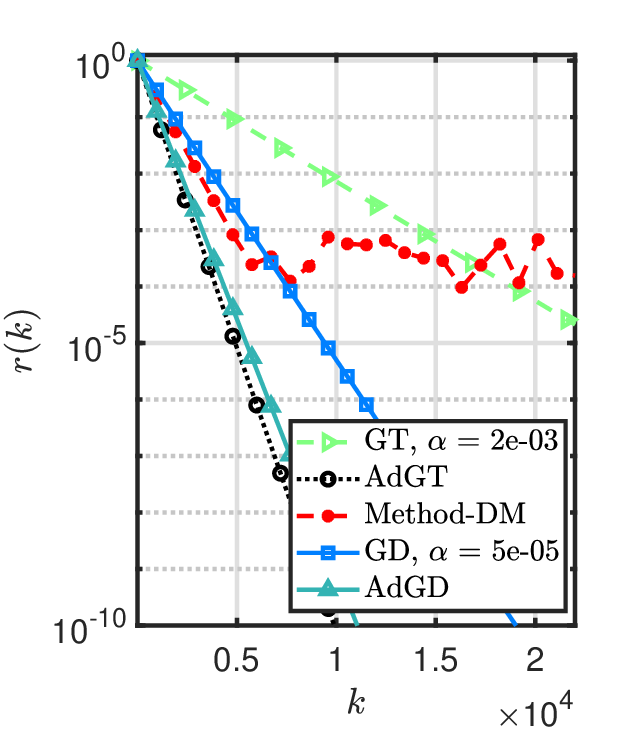}%
    }
    \hfill
    \subfloat[Scenario ii]{%
        \includegraphics[width=0.24\textwidth]{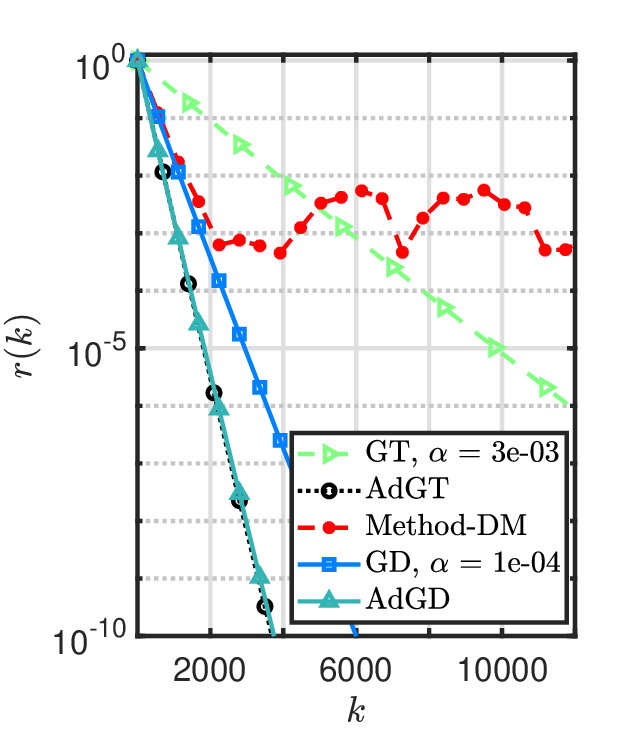}%
    }
    \hfill
    \subfloat[Scenario iii]{%
        \includegraphics[width=0.24\textwidth]{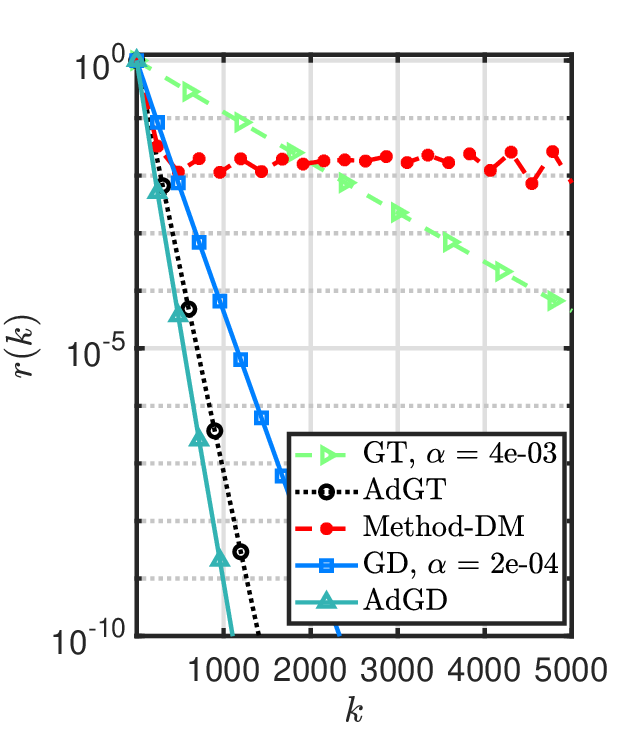}%
    }
    \hfill
    \subfloat[Scenario iv]{%
        \includegraphics[width=0.24\textwidth]{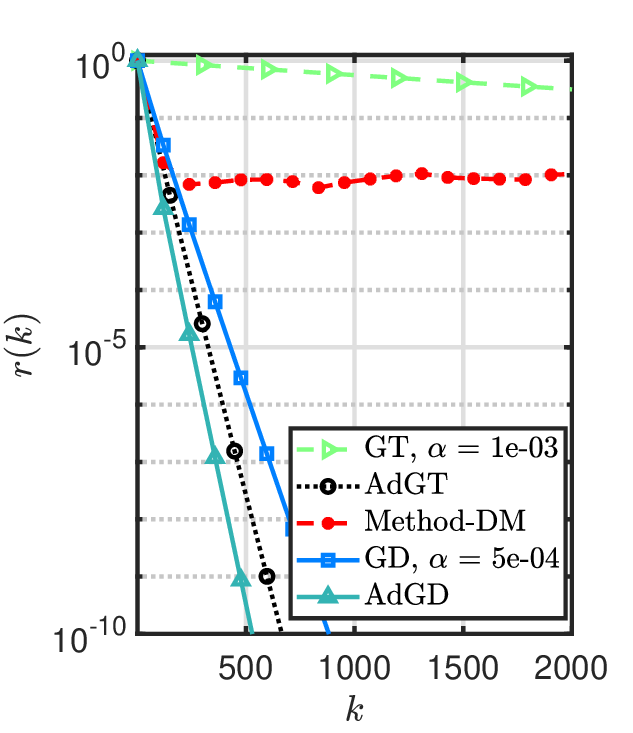}%
    }
    \caption{Decentralized quadratic problem with $n=100$ agents under four
heterogeneity scenarios.}
    \label{fig:quadratic}
\end{figure}
We {compare AdGT with GT and} Method-DM and {include centralized GD and AdGD as references. GD applies a fixed stepsize} to the aggregate objective {in~\eqref{Objective_quadratic}, whereas AdGD uses the adaptive update} of~\cite{malitsky2020adaptive}. {This comparison separates the general benefit of adaptive stepsizes from the additional benefit of per-agent adaptation in a decentralized network. We set $\gamma=1$ and select the} GD and GT {stepsizes for the best stable empirical performance over the tested grid. 

As} the number of {agents with $\tau=3$} decreases, the {separation between AdGT and tuned GT grows: AdGT reaches comparable residual levels in fewer iterations. In} Scenario~(iv){, the best stable GT stepsize found over the tested grid is approximately} $10^{-3}$, {whereas the AdGT agents with $\tau=1$} use stepsizes between $10^{-2}$ {and $10^{-1}$. By} contrast, the {AdGD--GD separation changes much less across the four scenarios. This indicates that per-agent adaptation reduces the bottleneck created when
one common decentralized stepsize must remain stable for a small number of
agents with large local smoothness constants}.

\begin{figure}[!ht]
    \centering
    \subfloat[Agent with $\tau=1$]{%
        \includegraphics[width=0.45\textwidth]{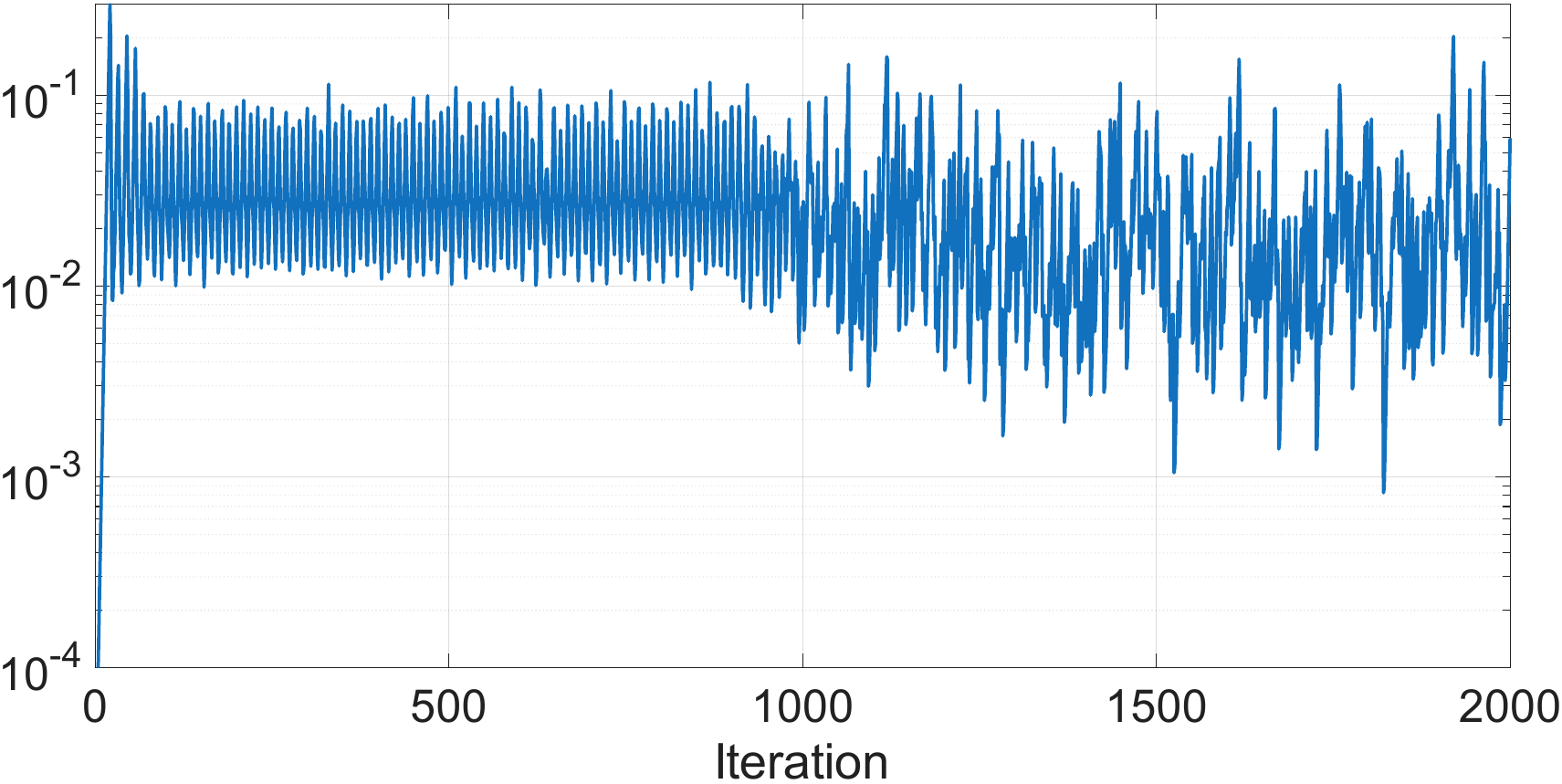}%
    }
    \hfill
    \subfloat[Agent with $\tau=3$]{%
        \includegraphics[width=0.45 \textwidth]{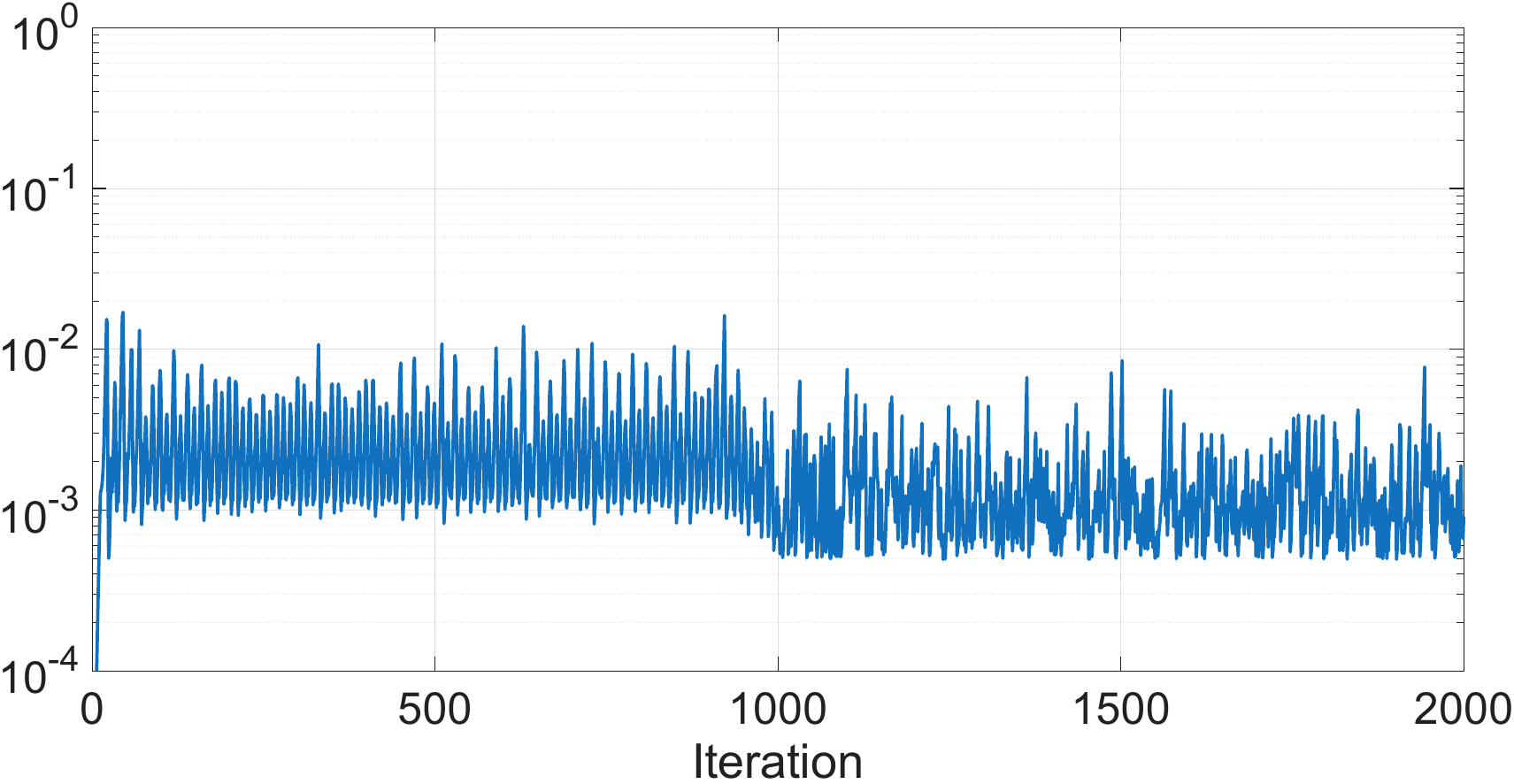}%
    }

\caption{{Stepsizes for two different agents in Scenario iv of Fig.~\ref{fig:quadratic}.}}
\label{Stepsise_agent_QF}
\end{figure}

\subsection{{Ridge Regression with Different Numbers} of {Agents} and {Network
Topologies}}\label{Simulation:section_3}
We consider ridge regression as a specific instance of problem~\eqref{eqz1}, where the local objective function for each agent~$i$ is defined as  
\begin{align*}
    f_i(\bm{x}) = \|A_i \bm{x} - \bm{b}_i\|^2 + \rho \|\bm{x}\|^2,
\end{align*}
with $A_i \in \mathbb{R}^{m \times p}$, $\bm{b}_i \in \mathbb{R}^{m}$, and the regularization parameter set to $\rho = 0.1$. The data used in this experiment corresponds to the \texttt{w8a} dataset, as described in Section~\ref{Simulation:section_1}.

To evaluate performance, we conducted numerical experiments on both random and cycle networks. For the random networks, we generated graphs with a connectivity ratio of $0.35$, considering two cases with $n = 25$ and $n = 100$ nodes. For the cycle graph, we tested it with $n = 10$ and $n = 25$ nodes.
\begin{figure}[!ht]
    \centering
    \subfloat[Cycle, $n=10$]{%
        \includegraphics[width=0.24\textwidth]{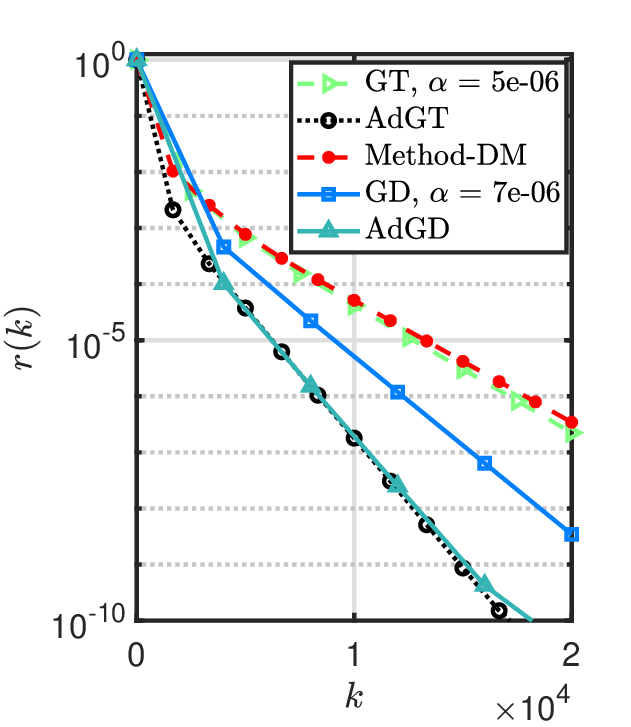}%
    }
    \hfill
    \subfloat[Cycle, $n=25$]{%
        \includegraphics[width=0.24\textwidth]{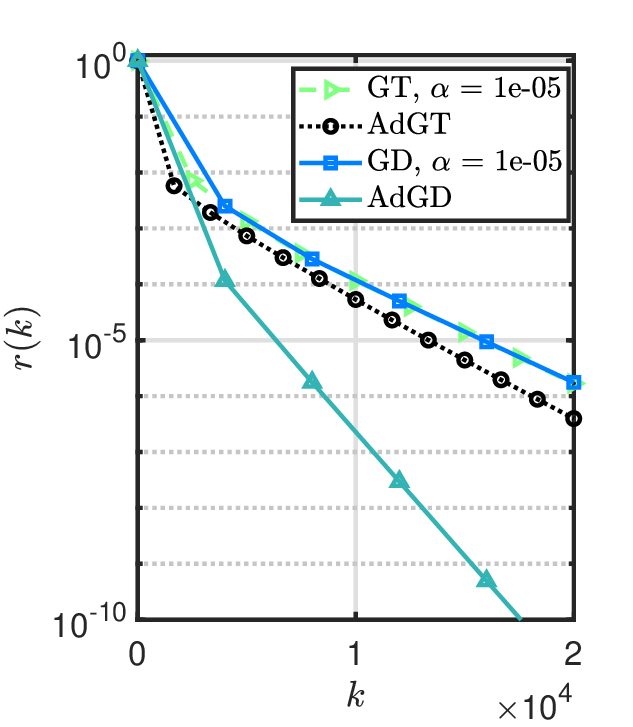}%
    }
    \hfill
    \subfloat[Random, $n=25$]{%
        \includegraphics[width=0.24\textwidth]{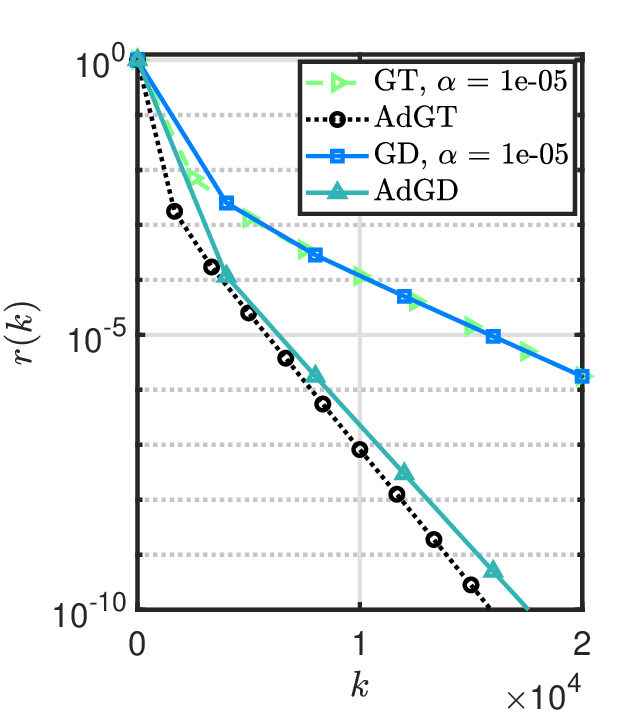}%
    }
    \hfill
    \subfloat[Random, $n=100$]{%
        \includegraphics[width=0.24\textwidth]{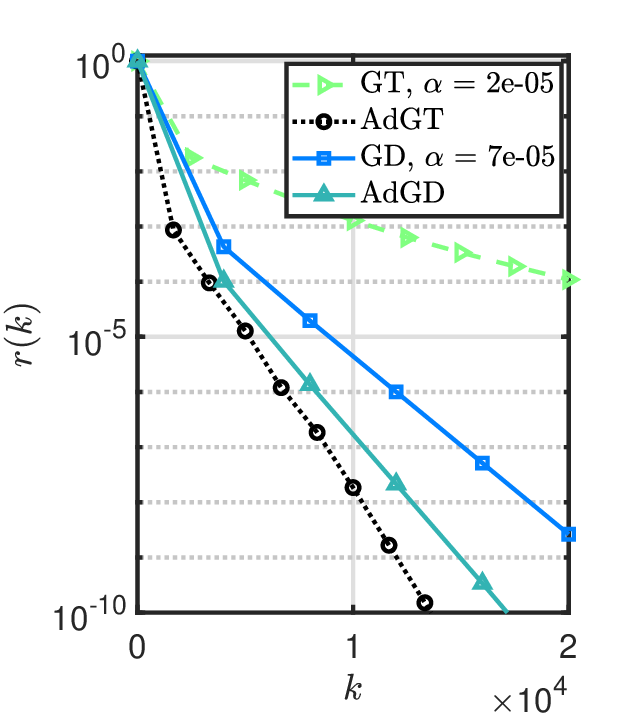}%
    }
    \caption{Ridge regression over cycle and random graphs with different numbers
of agents.}
    \label{Ridge_fig}
\end{figure}
In all four experimental settings, the full dataset was evenly partitioned among the $n$ agents. Thus, each agent receives approximately $m=\lfloor 49749/n \rfloor$ samples drawn without replacement from the global dataset. We note that these experiments vary both the number of agents and the per-agent sample size. Thus, the observed trends should not be interpreted as a pure connectivity effect; they reflect the combined influence of graph topology, number of agents, and local data partitioning. We fix $\gamma = 1$, and the stepsizes for GD and GT are carefully tuned to obtain their best stable convergence performance.

Figure~\ref{Ridge_fig} {reports the residuals. On} the random graph{, the AdGT--GT separation is larger for $n=100$ than for $n=25$, whereas on} the cycle graph {it is larger for $n=10$ than for $n=25$. Because changing $n$ also changes the local sample size, these results should not be interpreted as a pure connectivity effect.} Method-DM {reaches a residual comparable to AdGT only} on the cycle graph with {$n=10$. As} in Section~\ref{Simulation:section_2}, the difference in {the iterations needed by} AdGT and GT {to reach the same residual is larger than the corresponding AdGD--GD difference, especially on the random graph with $n=100$}.

\section{Conclusion} We proposed AdGT, a decentralized gradient-tracking method with adaptive per-agent stepsizes{. For update~\eqref{alpha_nound_L_F}, we proved exact linear convergence under} smoothness, strong convexity, and undirected-network assumptions{. For updates~\eqref{Stepsize_Adgt_y_i} and \eqref{Stepsize_Adgt_y_i_f_i}, we characterized when the tracking-based candidate determines the stepsize and established conditional lower and upper stepsize bounds and linear convergence under a relative tracking-disagreement condition. The experiments show reduced graph-dependent stepsize tuning. A key empirical finding is that the AdGT--GT separation grows when only a small subset of agents has large} local smoothness constants{, whereas the corresponding AdGD--GD separation changes much less. This suggests that per-agent adaptation is especially useful when a common decentralized stepsize is restricted by only a few agents}. Future work {will derive the tracking-disagreement condition from standard graph and objective assumptions and extend AdGT to directed,} time-varying{, stochastic, convex, and nonconvex settings.}

%\appendices
\appendix
\section{{Proofs} of Lemma~\ref{lemma:alpha_k_bounded} { and
Proposition~\ref{prop:tracking_bounds_convergence}}}
\label{app:proof_alpha_k_bounded}

\begin{proof}
Fix an agent $i\in\mathcal V$. At $k=0$, since $\theta_i^0=0$,
\[
\alpha_i^1
=
\min\left\{\frac{1}{2\gamma L_{f_i}^0},\;\alpha_i^0\right\}.
\]
By $L_i$-smoothness, $L_{f_i}^0\le L_i$, and by assumption
$\alpha_i^0\ge 1/(2\gamma L_i)$. Hence both candidates in the minimum are at
least $1/(2\gamma L_i)$, so
\[
\alpha_i^1\ge \frac{1}{2\gamma L_i}\ge \frac{1}{2\gamma L},
\]
where $L=\max_{j\in\mathcal V}L_j$.
Suppose now that $\alpha_i^k\ge 1/(2\gamma L_i)$ for some $k\ge1$. If
$\|\bm{x}_i^{k+1}-\bm{x}_i^k\|=0$, then $\alpha_i^{k+1}=\alpha_i^k$, and the
bound is preserved. Otherwise, $L_{f_i}^k$ is well defined and satisfies
$L_{f_i}^k\le L_i$. Moreover, $\theta_i^k\ge0$. Therefore both candidates
in~\eqref{alpha_nound_L_F} are at least $1/(2\gamma L_i)$, and hence
\[
\alpha_i^{k+1}\ge \frac{1}{2\gamma L_i}.
\]
By induction,
\[
\alpha_i^k\ge \frac{1}{2\gamma L_i}\ge \frac{1}{2\gamma L},
\qquad \forall k\ge1.
\]

We next prove the upper bound. Since $L_{f_i}^0$ is well defined and $f_i$ is
$\mu_i$-strongly convex, we have $L_{f_i}^0\ge\mu_i$. Thus
\[
\alpha_i^1
\le
\frac{1}{2\gamma L_{f_i}^0}
\le
\frac{1}{2\gamma\mu_i}
\le
\frac{1}{2\gamma\mu}.
\]
Suppose now that $\alpha_i^k\le 1/(2\gamma\mu)$ for some $k\ge1$. If
$\|\bm{x}_i^{k+1}-\bm{x}_i^k\|=0$, then $\alpha_i^{k+1}=\alpha_i^k$, and the
bound is preserved. Otherwise, strong convexity, together with Cauchy--Schwarz,
gives $L_{f_i}^k\ge\mu_i$, and from~\eqref{alpha_nound_L_F},
\[
\alpha_i^{k+1}
\le
\frac{1}{2\gamma L_{f_i}^k}
\le
\frac{1}{2\gamma\mu_i}
\le
\frac{1}{2\gamma\mu}.
\]
By induction,
\[
\alpha_i^k\le \frac{1}{2\gamma\mu},
\qquad \forall k\ge1.
\]
Combining the lower and upper bounds gives
\[
\frac{1}{2\gamma L}
\le
\alpha_i^k
\le
\frac{1}{2\gamma\mu},
\qquad \forall i\in\mathcal V,\ k\ge1.
\]
\end{proof}
{
\subsection{Proof of Proposition~\ref{prop:tracking_bounds_convergence}}
\label{app:proof_tracking_bounds}
\begin{proof}
Consider an iteration with
$\|\bm x_i^{k+1}-\bm x_i^k\|>0$. Using the notation of
Remark~\ref{rem:tracking_candidate}, the reverse triangle inequality and the
relative tracking-disagreement condition give
\[
 \big|L_{y_i}^k-L_{f_i}^k\big|
 \le
 \frac{\|\bm r_i^k\|}
      {\|\bm x_i^{k+1}-\bm x_i^k\|}
 \le
 \xi_i.
\]
Smoothness and strong convexity imply
\[
 \mu_i
 \le
 L_{f_i}^k
 \le
 L_i,
\]
and therefore
\[
 L_{y_i}^k
 \le
 L_i+\xi_i.
\]
If $\xi_i<\mu_i$, they also give
\[
 L_{y_i}^k
 \ge
 \mu_i-\xi_i
 >
 0.
\]

We first consider update~\eqref{Stepsize_Adgt_y_i_f_i}. The initialization
satisfies
\[
 \alpha_i^0
 \ge
 \frac{1}{2\gamma L_i}
 \ge
 \frac{1}{2\gamma(L_i+\xi_i)}.
\]
At every iteration with a nonzero local primal change,
\[
 \frac{1}{2\gamma L_{f_i}^k}
 \ge
 \frac{1}{2\gamma(L_i+\xi_i)},
 \qquad
 \frac{1}{2\gamma L_{y_i}^k}
 \ge
 \frac{1}{2\gamma(L_i+\xi_i)},
\]
and
\[
 \sqrt{1+\theta_i^k}\,\alpha_i^k
 \ge
 \alpha_i^k.
\]
Thus, induction and the zero-movement convention give
\[
 \alpha_i^k
 \ge
 \frac{1}{2\gamma(L_i+\xi_i)},
 \qquad k\ge1.
\]
The retained local-gradient candidate also gives
\[
 \alpha_i^{k+1}
 \le
 \frac{1}{2\gamma L_{f_i}^k}
 \le
 \frac{1}{2\gamma\mu_i}.
\]
Hence, for update~\eqref{Stepsize_Adgt_y_i_f_i},
\[
 \frac{1}{2\gamma(L_i+\xi_i)}
 \le
 \alpha_i^k
 \le
 \frac{1}{2\gamma\mu_i}.
\]

Taking common bounds over all agents gives
\[
 \underline\alpha
 =
 \frac{1}{2\gamma\max_i(L_i+\xi_i)}
 \le
 \alpha_i^k
 \le
 \frac{1}{2\gamma\mu}.
\]
The consensus and gradient-tracking recursions in
Lemma~\ref{lemma:one_step_recursion} use the common upper stepsize bound, while
the optimality recursion also uses the positive lower bound. Replacing these
two constants therefore gives the comparison matrix stated in the proposition.
If $\gamma\ge L/\mu$, its upper stepsize is at most $1/(2L)$. When the modified
matrix has spectral radius below one, the norm argument of
Theorem~\ref{theorem2} gives linear convergence with a rate arbitrarily close
to that spectral radius.

We next consider update~\eqref{Stepsize_Adgt_y_i}. Its lower-bound proof is
unchanged because the tracking-based and growth-control candidates satisfy the
same lower bounds used above. The difference is the upper bound. If
$\xi_i<\mu_i$, then
\[
 \alpha_i^{k+1}
 \le
 \frac{1}{2\gamma L_{y_i}^k}
 \le
 \frac{1}{2\gamma(\mu_i-\xi_i)}.
\]
Consequently,
\[
 \frac{1}{2\gamma(L_i+\xi_i)}
 \le
 \alpha_i^k
 \le
 \frac{1}{2\gamma(\mu_i-\xi_i)}.
\]
The common upper bound is now
\[
 \frac{1}{2\gamma\min_i(\mu_i-\xi_i)}.
\]
Replacing only the upper stepsize bound in the preceding comparison argument
gives the matrix stated for update~\eqref{Stepsize_Adgt_y_i}. The condition
\[
 \gamma
 \ge
 \frac{L}{\min_i(\mu_i-\xi_i)}
\]
makes this upper bound at most $1/(2L)$, and the same spectral-radius argument
then gives linear convergence.
\end{proof}}

\section{Auxiliary Bounds for Lemma~\ref{lemma:one_step_recursion}}
\label{app:proof_one_step_recursion}

In this appendix, we provide the detailed derivation of the three inequalities summarized in Lemma~\ref{lemma:one_step_recursion}. We first establish an auxiliary bound on $\|\bm{\mathrm{y}}^k\|$, which is used in the consensus and tracking recursions.
Because $\bm{\mathrm{y}}^0=\nabla f(\bm{\mathrm{x}}^0)$ and $W$ is doubly stochastic, the
gradient-tracking recursion preserves the average-gradient identity
\begin{equation}\label{eq:appendix_ybar_identity}
    \bar{\bm{\mathrm{y}}}^k
    =
    \frac{1}{n}\bm{1}\bm{1}^\top \nabla f(\bm{\mathrm{x}}^k),
    \qquad \forall k\ge 0.
\end{equation}
Indeed, multiplying~\eqref{eqComp2} from the left by $\frac{1}{n}\bm{1}^\top$ gives
\[
\bar{\bm{y}}^{k+1}
=
\bar{\bm{y}}^k
+
\frac{1}{n}\bm{1}^\top
\bigl(\nabla f(\bm{\mathrm{x}}^{k+1})-\nabla f(\bm{\mathrm{x}}^k)\bigr),
\]
where $\bar{\bm{y}}^k=\frac{1}{n}\bm{1}^\top\bm{\mathrm{y}}^k$. Since
$\bar{\bm{y}}^0=\frac{1}{n}\bm{1}^\top\nabla f(\bm{\mathrm{x}}^0)$, the claim follows by
induction and by multiplying the row identity by $\bm{1}$.
%\subsection{Bound on $\|\bm{\mathrm{y}}^k\|$}

\subsection{Bound on \texorpdfstring{$\|\bm{\mathrm{y}}^k\|$}{norm of y to the k}}
\label{app:proof_bound_y}

For all $k$, the following inequality holds:
\begin{align}
    \|\bm{\mathrm{y}}^k\|\leq
         \|\bm{\mathrm{y}}^k- \bar{\bm{\mathrm{y}}}^k \| 
         +  L\|\bm{\mathrm{x}}^k- \bar{\bm{\mathrm{x}}}^k\|  
         +\sqrt{n} L\| \bar{\bm{{x}}}^k - \bm{{x}}^* \|.
\label{eq:appendix_y_bound}
\end{align}

\begin{proof}
We add and subtract $\tfrac{1}{n}\bm{1}\bm{1}^\top \nabla f(\bm{\mathrm{x}}^k)$, use the fact that $\tfrac{1}{n}\bm{1}\bm{1}^\top  \nabla f(\bm{x}^*)=0$, and apply the triangle inequality:
\begin{align*}
    \|\bm{\mathrm{y}}^k\|
    &\leq 
     \|\bm{\mathrm{y}}^k-\tfrac{1}{n} \bm{1}\bm{1}^\top \nabla f(\bm{\mathrm{x}}^k)\| 
     + \|\tfrac{1}{n} \bm{1}\bm{1}^\top \nabla f(\bm{\mathrm{x}}^k)-\tfrac{1}{n} \bm{1}\bm{1}^\top \nabla f(\bar{\bm{x}}^{k})\|  \\
    &\quad +\|\tfrac{1}{n} \bm{1}\bm{1}^\top \nabla f(\bar{\bm{x}}^{k})-\tfrac{1}{n} \bm{1}\bm{1}^\top \nabla f(\bm{x}^*)\|.
\end{align*}
Using~\eqref{eq:appendix_ybar_identity} and $\|\bm{1}\bm{b}\|=\sqrt{n}\|\bm{b}\|$, this becomes
\begin{align*}
    \|\bm{\mathrm{y}}^k\|
    &\leq \|\bm{\mathrm{y}}^k- \bar{\bm{\mathrm{y}}}^k \|
     + \sqrt{n}\|\tfrac{1}{n}\bm{1}^\top \nabla f(\bm{\mathrm{x}}^k)-\tfrac{1}{n}\bm{1}^\top \nabla f(\bar{\bm{x}}^{k})\| \\
    &\quad + \sqrt{n}\|\tfrac{1}{n}\bm{1}^\top \nabla f(\bar{\bm{x}}^{k})-\tfrac{1}{n}\bm{1}^\top \nabla f(\bm{x}^*)\|.
\end{align*}
By $L$-smoothness,
\begin{align*}
    \|\bm{\mathrm{y}}^k\|
    \leq \|\bm{\mathrm{y}}^k- \bar{\bm{\mathrm{y}}}^k \|
    + \sqrt{n} L \sum_{i=1}^n \tfrac{\|\bm{{x}}_i^k-\bar{\bm{{x}}}^k\|}{n}
    + \sqrt{n} L \|\bar{\bm{{x}}}^k-\bm{{x}}^*\|.
\end{align*}
Finally, Jensen’s inequality gives
{
\[
\sum_{i=1}^n \tfrac{\|\bm{{x}}_i^k-\bar{\bm{{x}}}^k\|}{n}
\leq \sqrt{\sum_{i=1}^n \tfrac{\|\bm{{x}}_i^k-\bar{\bm{{x}}}^k\|^2}{n}}
= \tfrac{1}{\sqrt{n}}\|\bm{\mathrm{x}}^k-\bar{\bm{\mathrm{x}}}^k\|,
\]
}
which completes the proof.
\end{proof}

\subsection{Proof of the Consensus Recursion}
\label{app:proof_bound_x_diff}

We prove that, for all $k\geq1$,
{
\begin{align}
    \|\bm{\mathrm{x}}^{k+1}- \bar{\bm{\mathrm{x}}}^{k+1}\| 
    &\leq (\lambda+\lambda \alpha_{\max} L)\|\bm{\mathrm{x}}^k- \bar{\bm{\mathrm{x}}}^k\|  
    + \lambda \alpha_{\max} \|\bm{\mathrm{y}}^k- \bar{\bm{\mathrm{y}}}^k \|  + \sqrt{n}\,\lambda \alpha_{\max} L\| \bar{\bm{{x}}}^k - \bm{{x}}^* \|.
\label{eq:appendix_consensus_recursion}
\end{align}
}

\begin{proof}
From the update rule we have
\begin{align*}
    \bm{\mathrm{x}}^{k+1}- \bar{\bm{\mathrm{x}}}^{k+1}
    &= \left(W- \tfrac{1}{n}\bm{1}\bm{1}^\top\right)\big(\bm{\mathrm{x}}^k - D^k \bm{\mathrm{y}}^k\big).
\end{align*}
Taking norms and applying the triangle inequality, and using the fact that
$\widetilde{W}\bar{\bm{\mathrm{x}}}^k = 0$, yields
\begin{align*}
    \|\bm{\mathrm{x}}^{k+1}- \bar{\bm{\mathrm{x}}}^{k+1}\|
    &\leq \|\tilde{W}(\bm{\mathrm{x}}^k - \bar{\bm{\mathrm{x}}}^k)\|
     + \|\tilde{W}(D^k\bm{\mathrm{y}}^k)\|,
\end{align*}
where $\tilde{W} := W-\tfrac{1}{n}\bm{1}\bm{1}^\top$.  
Since $\vertiii{\tilde{W}}\leq \lambda$ and $D^k$ is diagonal with entries $\alpha_i^k \leq \alpha_{\max}$, we have $\vertiii{D^k}\leq \alpha_{\max}$, and hence
\begin{align*}
    \|\bm{\mathrm{x}}^{k+1}- \bar{\bm{\mathrm{x}}}^{k+1}\|
    &\leq \lambda \|\bm{\mathrm{x}}^k - \bar{\bm{\mathrm{x}}}^k\|
     + \lambda \alpha_{\max}\|\bm{\mathrm{y}}^k\|.
\end{align*}
Finally, applying \eqref{eq:appendix_y_bound} to bound $\|\bm{\mathrm{y}}^k\|$ yields the result.
\end{proof}

\subsection{Proof of the Gradient-Tracking Recursion}
\label{app:proof_bound_y_diff}

We prove that, for all $k\geq1$,
{
\begin{align}
   \| \bm{\mathrm{y}}^{k+1}- \bar{\bm{\mathrm{y}}}^{k+1} \| 
   &\leq (\lambda+L\alpha_{\max}) \|\bm{\mathrm{y}}^k- \bar{\bm{\mathrm{y}}}^k \| 
    +  L \big(\vertiii{W-I}+ L\alpha_{\max}\big)\|\bm{\mathrm{x}}^k- \bar{\bm{\mathrm{x}}}^k\|  \nonumber\\
   &\quad +\sqrt{n}\, L^2 \alpha_{\max} \| \bar{\bm{x}}^k - \bm{x}^* \|.
\label{eq:appendix_tracking_recursion}
\end{align}
}

\begin{proof}
From the update rule we can write
\begin{align*}
    \bm{\mathrm{y}}^{k+1}- \bar{\bm{\mathrm{y}}}^{k+1} 
    &=   W \bm{\mathrm{y}}^k 
        +  \nabla f(\bm{\mathrm{x}}^{k+1})-\nabla f(\bm{\mathrm{x}}^k)  - \bar{\bm{\mathrm{y}}}^k 
        - \tfrac{1}{n} \bm{1}\bm{1}^\top \big(\nabla f(\bm{\mathrm{x}}^{k+1})-\nabla f(\bm{\mathrm{x}}^k)\big).
\end{align*} 
Taking norms and applying the triangle inequality gives
{
\begin{align}
     \|\bm{\mathrm{y}}^{k+1}- \bar{\bm{\mathrm{y}}}^{k+1} \| 
     &\leq    \|\tilde{W} (\bm{\mathrm{y}}^k - \bar{\bm{\mathrm{y}}}^{k}) \|  +  \big\|(I-\tfrac{1}{n} \bm{1}\bm{1}^\top)(\nabla f(\bm{\mathrm{x}}^{k+1})-\nabla f(\bm{\mathrm{x}}^k))\big\| \nonumber\\
     &\leq \lambda \|\bm{\mathrm{y}}^k - \bar{\bm{\mathrm{y}}}^{k}\|
     + L \|\bm{\mathrm{x}}^{k+1}-\bm{\mathrm{x}}^k\|, \label{eq:bound_diff_x_y_app}
\end{align}
}
where we used that the spectral norm of $\tilde{W}$ satisfies $\vertiii{\tilde{W}}\leq \lambda$, the spectral norm of $\vertiii{I-\tfrac{1}{n}\bm{1}\bm{1}^\top}$ equals $1$, and $L$-smoothness.  
Next, we bound the term $\|\bm{\mathrm{x}}^{k+1}-\bm{\mathrm{x}}^k\|$:
\begin{align}
  \|\bm{\mathrm{x}}^{k+1}-\bm{\mathrm{x}}^k\|
  &= \|(W-I)\bm{\mathrm{x}}^{k} - W D^k \bm{\mathrm{y}}^k\| \nonumber\\
  &\leq \|(W-I)(\bm{\mathrm{x}}^{k}-\bar{\bm{\mathrm{x}}}^{k})\| + \|W D^k \bm{\mathrm{y}}^k\| \nonumber\\
  &\leq \vertiii{W-I}\,\|\bm{\mathrm{x}}^{k}-\bar{\bm{\mathrm{x}}}^{k}\|
       + \alpha_{\max}\|\bm{\mathrm{y}}^k\|, \label{eq:diff_bound_x_k_app}
\end{align}
where in the second equation we used $(W-I)\bar{\bm{\mathrm{x}}}^k=0$, and in the last inequality we used that the spectral norm of $W$ equals $1$, together with the fact that $D^k$ is diagonal with entries $\alpha_i^k \leq \alpha_{\max}$, hence $\vertiii{D^k}\leq \alpha_{\max}$.\footnote{In particular, if $W$ is doubly stochastic, then $\vertiii{W-I}\leq 2$.}
Finally, substituting \eqref{eq:diff_bound_x_k_app} into \eqref{eq:bound_diff_x_y_app}, and applying \eqref{eq:appendix_y_bound} to bound $\|\bm{\mathrm{y}}^k\|$, yields the stated inequality.
\end{proof}

\subsection{Proof of the Averaged-Iterate Recursion}
\label{app:proof_barx_contraction}

We prove that if $\gamma \geq \tfrac{L}{\mu}$, then, for all $k\geq1$,
{
\begin{align}
\| \bar{\bm{x}}^{k+1} - \bm{x}^*\|  &\leq  \left(1- \frac{\mu}{2 \gamma L}\right) \,\|\bar{\bm{x}}^{k} -\bm{x}^{*}\| 
     + L \frac{\alpha_{\max}}{\sqrt{n}}  \|\bar{\bm{\mathrm{x}}}^{k} -\bm{\mathrm{x}}^{k}  \|  +\frac{\alpha_{\max}}{\sqrt{n}} \|\bm{\mathrm{y}}^{k} -\bar{\bm{\mathrm{y}}}^{k}\| .
\label{eq:appendix_average_recursion}
\end{align}
}

\begin{proof}
Taking the average of \eqref{eqComp1}, we obtain
\[
\bar{\bm{x}}^{k+1}-\bm{x}^*
= \bar{\bm{x}}^{k}-\tfrac{1}{n}\bm{1}^\top D^k \bm{\mathrm{y}}^{k}-\bm{x}^*.
\]
Adding and subtracting $\tfrac{1}{n}\bm{1}^\top D^k\bar{\bm{\mathrm{y}}}^{k}$ and applying the triangle inequality yield
\begin{align}
      \|\bar{\bm{x}}^{k+1} - \bm{x}^*\| &\leq  \|\bar{\bm{x}}^{k} - \tfrac{1}{n} \bm{1}^\top D^k \bar{\bm{\mathrm{y}}}^{k} -  \bm{x}^* \| +\| \tfrac{1}{n} \bm{1}^\top D^k (\bm{\mathrm{y}}^{k} -\bar{\bm{\mathrm{y}}}^{k})
      \| \nonumber\\
      &\leq  \|\bar{\bm{x}}^{k} - \tfrac{1}{n} \bm{1}^\top D^k \bar{\bm{\mathrm{y}}}^{k} -  \bm{x}^* \|  +\| \tfrac{1}{n} \bm{1}^\top D^k\| \| \bm{\mathrm{y}}^{k} -\bar{\bm{\mathrm{y}}}^{k}
      \| \nonumber\\
       &\leq  \|\bar{\bm{x}}^{k} - \tfrac{1}{n} \bm{1}^\top D^k \bar{\bm{\mathrm{y}}}^{k} -  \bm{x}^* \|  +\frac{\alpha_{\max}}{\sqrt{n}} \|\bm{\mathrm{y}}^{k} -\bar{\bm{\mathrm{y}}}^{k}\| . \label{eq:barx_step_split_app}
  \end{align}
In the last inequality we used the spectral norm $\|\tfrac{1}{n}\bm{1}^\top D^k\|\le \alpha_{\max}/\sqrt{n}$.\\
Next, add and subtract $\tfrac{1}{n}\bm{1}^\top D^k\,\tfrac{1}{n}\bm{1}\bm{1}^\top\nabla f(\bar{\bm{x}}^{k})$ and use the identity~\eqref{eq:appendix_ybar_identity}:
 \begin{align}
     & \|\bar{\bm{x}}^{k} - \tfrac{1}{n} \bm{1}^\top D^k \bar{\bm{\mathrm{y}}}^{k} -  \bm{x}^* \| \nonumber\\
      &\leq \|\bar{\bm{x}}^{k} - \tfrac{1}{n} \bm{1}^\top D^k \tfrac{1}{n}\bm{1}\bm{1}^\top \nabla f(\bar{\bm{x}}^{k}) -  \bm{x}^* \| +\| \tfrac{1}{n} \bm{1}^\top D^k \tfrac{1}{n}\bm{1}\bm{1}^\top \nabla f(\bar{\bm{x}}^{k}) -\tfrac{1}{n} \bm{1}^\top D^k  \bar{\bm{\mathrm{y}}}^{k} \| \nonumber\\
      &\leq \|\bar{\bm{x}}^{k} - \big(\tfrac{1}{n} \sum_{i=1}^n \alpha_i^k \big) \big(\tfrac{1}{n} \sum_{i=1}^n \nabla f_i(\bar{\bm{x}}^{k}) \big) -  \bm{x}^* \|  +\| \tfrac{1}{n} \bm{1}^\top D^k \tfrac{1}{n}\bm{1}\bm{1}^\top(\nabla f(\bar{\bm{x}}^{k}) - \nabla f(\bm{\mathrm{x}}^{k}) ) \| \nonumber\\
      &\leq \|\bar{\bm{x}}^{k} - \big(\tfrac{1}{n} \sum_{i=1}^n \alpha_i^k \big)\big(\tfrac{1}{n} \sum_{i=1}^n \nabla f_i(\bar{\bm{x}}^{k}) \big) -  \bm{x}^* \|  +L \frac{\alpha_{\max}}{\sqrt{n}}  \|\bar{\bm{\mathrm{x}}}^{k} -\bm{\mathrm{x}}^{k}  \| .\label{eq:centralized_plus_consensus_app}
  \end{align}
The last line follows from $L$-smoothness and the same bound on $\tfrac{1}{n}\bm{1}^\top D^k$.
When $\tfrac{1}{n}\sum_{i=1}^n \alpha_i^k<\tfrac{2}{L}$, Lemma~\ref{lemma_max_L_mu} gives
\[
\Bigl\|\bar{\bm{x}}^{k}-\Bigl(\tfrac{1}{n}\sum_{i=1}^n \alpha_i^k\Bigr)\Bigl(\tfrac{1}{n}\sum_{i=1}^n \nabla f_i(\bar{\bm{x}}^{k})\Bigr)-\bm{x}^*\Bigr\|
\le \eta\,\|\bar{\bm{x}}^{k}-\bm{x}^*\|,
\]
where
{
\[
\eta \triangleq \max\!\left\{\,\bigl|1-L\,\tfrac{1}{n}\!\sum_{i=1}^n \alpha_i^k\bigr|,\;
\bigl|1-\mu\,\tfrac{1}{n}\!\sum_{i=1}^n \alpha_i^k\bigr|\,\right\}.
\]
}
By Lemma~\ref{lemma:alpha_k_bounded}, this stepsize condition holds whenever 
$\gamma>\tfrac{L}{4 \mu}$. Moreover, under $\gamma\ge \tfrac{L}{\mu}$, both terms inside the maximum
are positive and strictly less than $1$, so the maximum reduces to
\[
\eta=1-\mu\Bigl(\tfrac{1}{n}\sum_{i=1}^n \alpha_i^k\Bigr).
\]
Finally, since Lemma~\ref{lemma:alpha_k_bounded} guarantees $\alpha_i^k\ge \tfrac{1}{2\gamma L}$, we obtain
\begin{equation}\label{eq:effect_1}
 \eta \le 1-\frac{\mu}{2\gamma L}.   
\end{equation}
Combining \eqref{eq:centralized_plus_consensus_app} with the bound on $\eta$ yields
\begin{align}
\Bigl\|\bar{\bm{x}}^{k}-\tfrac{1}{n}\bm{1}^\top D^k \bar{\bm{\mathrm{y}}}^{k}-\bm{x}^*\Bigr\|
&\le \Bigl(1-\frac{\mu}{2\gamma L}\Bigr)\|\bar{\bm{x}}^{k}-\bm{x}^*\|  + L\,\frac{\alpha_{\max}}{\sqrt{n}}\;\|\bar{\bm{\mathrm{x}}}^{k}-\bm{\mathrm{x}}^{k}\|.
\label{eq:barx_first_term_final_app}
\end{align}
Finally, substituting \eqref{eq:barx_first_term_final_app} into \eqref{eq:barx_step_split_app} gives
\begin{align*}
\|\bar{\bm{x}}^{k+1}-\bm{x}^*\|
& \leq \Bigl(1-\frac{\mu}{2\gamma L}\Bigr)\|\bar{\bm{x}}^{k}-\bm{x}^*\|
   + L\,\frac{\alpha_{\max}}{\sqrt{n}}\;\|\bar{\bm{\mathrm{x}}}^{k}-\bm{\mathrm{x}}^{k}\| + \frac{\alpha_{\max}}{\sqrt{n}}\;\|\bm{\mathrm{y}}^{k}-\bar{\bm{\mathrm{y}}}^{k}\|,
\end{align*}
which proves the claim.
\end{proof}

Combining \eqref{eq:appendix_consensus_recursion}, \eqref{eq:appendix_tracking_recursion}, and \eqref{eq:appendix_average_recursion} proves Lemma~\ref{lemma:one_step_recursion}.

\section{Proof of Theorem~\ref{theorem1}}
\label{app:proof_theorem1}

\begin{proof}
Given $\gamma \geq \frac{L}{\mu}$, it follows from Lemma~\ref{lemma:one_step_recursion} that \eqref{eqz19} holds for $k\geq1$. Next, we show $\rho(\Upsilon)<1$.
Since $\Upsilon$ has all non-negative entries, 
it is sufficient to show that $\Upsilon~ \beta< \beta$ for some positive $\beta=\left[\beta_{1},\beta_{2},\beta_{3} \right]^{\top} \in \mathbb{R}^{3}_+$; see \cite[Corollary~8.1.29]{horn2012matrix}. 
Hence, $\Upsilon~\beta< \beta$ is equivalent to
{
\begin{subequations}\label{eqqz22}
\begin{align}
(\lambda L\beta_{1}+\lambda\beta_{2}+\sqrt{n}\,\lambda L\beta_{3})\, \alpha_{\max} &< \beta_{1}(1-\lambda), \label{eqqz22a}\\
(L^2\beta_{1}+ L\beta_{2}+ \sqrt{n}\, L^2 \beta_{3})\,\alpha_{\max}  &< \beta_{2}(1-\lambda)-L\vertiii{W-I}\, \beta_1 , \label{eqqz22b}\\
\left(\frac{L}{\sqrt{n}}\beta_{1}+\frac{1}{\sqrt{n}}\beta_{2} \right)\alpha_{\max} &< \frac{\mu \beta_{3}}{2 \gamma L}. \label{eqqz22c}
\end{align}
\end{subequations}
}
For \eqref{eqqz22b} to be feasible, the right-hand side must be positive, which enforces
{
\begin{equation}\label{eq:beta_1_bound}
   0 \;<\; \beta_1 \;<\; \frac{(1-\lambda)}{L\vertiii{W-I}}\,\beta_2 .
\end{equation}
}
Moreover, by \eqref{define_bound_alpha_max} we have the upper bound 
\(\alpha_{\max} \le \frac{1}{2\gamma\mu}\). Substituting this value into
\eqref{eqqz22c} yields
\begin{equation}\label{eq:beta_3_bound}
   \beta_3 \;>\; \left(\frac{L}{\sqrt{n}}\beta_{1}+\frac{1}{\sqrt{n}}\beta_{2}\right)\frac{L}{\mu^{2}}
\end{equation}
is sufficient. Let $\delta_1 = \tfrac{\beta_1}{\beta_2}$ and $\delta_2 = \tfrac{\beta_3}{\beta_2}$ for any
fixed $\beta_2>0$. Under \eqref{eq:beta_1_bound} and \eqref{eq:beta_3_bound}, the inequalities \eqref{eqqz22a}--\eqref{eqqz22b} are guaranteed whenever
\begin{equation}\label{eq:alpha_max_bar_alpha_b}
   \alpha_{\max} \;<\; \bar{\alpha},
\end{equation}
where $\bar\alpha$ is defined in \eqref{eq:alpha_bar_def}. Because the upper bound 
\(\alpha_{\max} \le \frac{1}{2\gamma\mu}\), condition \eqref{eq:alpha_max_bar_alpha_b} follows from $\gamma>1/(2\mu\bar\alpha)$. Combining this with the standing requirement $\gamma>L/\mu$ ensures the existence of $\beta\in\mathbb{R}^3_{++}$ with $\Upsilon\beta<\beta$, hence $\rho(\Upsilon)<1$. 
\end{proof}

\section{Proof of Theorem~\ref{theorem2}}
\label{app:proof_theorem2}

\begin{proof}
By Theorem~\ref{theorem1}, $\rho(\Upsilon)<1$. Hence, by
\cite[Lemma~5.6.10]{horn2012matrix}, for any
$\zeta\in(0,\,1-\rho(\Upsilon))$ there exists an induced matrix norm
$\|\cdot\|_{(\zeta)}$ such that
\[
    \|\Upsilon\|_{(\zeta)}\le \rho(\Upsilon)+\zeta<1 .
\]
Since the comparison recursion~\eqref{eqz19} holds for all $k\geq1$, iterating it gives
\[
    e_k \le \Upsilon^{k-1}e_1, \qquad \forall k\geq1,
\]
where the inequality is componentwise. Since all entries of $e_k$ and
$\Upsilon^{k-1}e_1$ are nonnegative and the Euclidean norm is monotone,
componentwise domination implies
\[
    \|e_k\| \le \|\Upsilon^{k-1}e_1\|.
\]
By equivalence of norms in finite dimension, there exists a constant
$\Gamma>0$ such that
\begin{align*}
\|e_k\|
\le
\Gamma\,\|\Upsilon^{k-1}e_1\|_{(\zeta)}
\le&
\Gamma\,\|\Upsilon\|_{(\zeta)}^{k-1}\|e_1\|_{(\zeta)}\\
\le&
\Gamma\,\bigl(\rho(\Upsilon)+\zeta\bigr)^{k-1}\|e_1\|,
\qquad \forall k\geq1,
\end{align*}

where $\Gamma$ absorbs the norm-equivalence constants.
Finally, for all $k\geq1$,
\[
\|\bm{\mathrm{x}}^k-\bm{1}\bm{x}^*\|
\le
\|\bm{\mathrm{x}}^k-\bar{\bm{\mathrm{x}}}^k\|
+
\sqrt{n}\,\|\bar{\bm{x}}^k-\bm{x}^*\|
\le
(1+\sqrt{n})\,\|e_k\|.
\]
Therefore, for all $k\geq1$,
\[
\|\bm{\mathrm{x}}^k-\bm{1}\bm{x}^*\|
\le
(1+\sqrt{n})\,\Gamma\,\|e_1\|
\bigl(\rho(\Upsilon)+\zeta\bigr)^{k-1}.
\]
Since $\zeta\in(0,1-\rho(\Upsilon))$ is arbitrary, the sequence
$\{\bm{\mathrm{x}}^k\}$ generated by AdGT converges exactly to
$\bm{\mathrm{x}}^*=\bm{1}\bm{x}^*$ after the initialization step, with a
geometric rate arbitrarily close to $\rho(\Upsilon)$.
\end{proof}

\bibliographystyle{IEEEtran}
\bibliography{references}
\end{document}